\documentclass[11pt, psamsfonts]{amsart}
\usepackage{amssymb}
\usepackage{array,calc,url,tabularx, epic,tikz}
\usetikzlibrary{arrows,positioning,decorations.pathmorphing,
  decorations.markings}

\def\ignore#1{\relax}

\def\g{\mathfrak g}

\def\R{{\mathbb R}}
\def\Z{{\mathbb Z}}

\def\C{{\mathbb C}}

\def\la{\lambda}
\def\La{\Lambda}

\def\N{\mathbb N}

\def\Ca{\mathcal C}

\def\P{\mathcal P}
\def\Pn{{\mathcal P}_n}

\def\B{{\mathcal B}}

\def\U{{\bf U}}

\def\ignore#1{\relax}

\def\om{\omega}
\def\eps{\epsilon}
\def\1{{\bf 1}}

\def\End{{\rm End}}
\def\Hom{{\rm Hom}}

\def\eps{\varepsilon}

\calclayout

\renewenvironment{proof}[1][\proofname]{\par
  \normalfont
  \topsep6\p@\@plus6\p@ \trivlist
  \item[\hskip\labelsep\bfseries
    #1\@addpunct{.}]\ignorespaces
}{%
  \qed\endtrivlist
}

\def\th@plain{%
  \let\thmhead\thmhead@plain \let\swappedhead\swappedhead@plain
  \thm@preskip.5\baselineskip\@plus.2\baselineskip
                                    \@minus.2\baselineskip
  \thm@postskip\thm@preskip
  \itshape
\renewcommand{\labelenumi}{{(\alph{enumi})\quad}}
                        \renewcommand{\labelenumii}{{(\roman{enumii})\ }}
}
\def\th@definition{%
  \let\thmhead\thmhead@plain \let\swappedhead\swappedhead@plain
  \thm@preskip.5\baselineskip\@plus.2\baselineskip
                                    \@minus.2\baselineskip
  \thm@postskip\thm@preskip
  \upshape
}
\def\th@remark{%
  \thm@headfont{\itshape}
  \let\thmhead\thmhead@plain \let\swappedhead\swappedhead@plain
  \thm@preskip.5\baselineskip\@plus.2\baselineskip
                                    \@minus.2\baselineskip
  \thm@postskip\thm@preskip
  \upshape
}

{\theoremstyle{plain}
\newtheorem{theorem}{Theorem}[section]
}

{\theoremstyle{plain}
\newtheorem{proposition}[theorem]{Proposition}
}

{\theoremstyle{plain}
\newtheorem{corollary}[theorem]{Corollary}
}

{\theoremstyle{plain}
\newtheorem{lemma}[theorem]{Lemma}
}

{\theoremstyle{plain}

}

{\theoremstyle{definition}
\newtheorem{definition}[theorem]{Definition}
}

{\theoremstyle{definition}
\newtheorem{example}[theorem]{Example}
}

{\theoremstyle{remark}
\newtheorem{remark}[theorem]{Remark}
}

{\theoremstyle{remark}

}

\numberwithin{equation}{section}

\renewcommand{\labelenumi}{{ \theenumi.}}
\renewcommand{\labelenumii}{{(\alph{enumii})}}

\def\la{\lambda}
\def\al{\alpha}

\def\UU{{\mathcal U}}

\def\choose #1 #2{\begin{pmatrix}#1\\#2\end{pmatrix}}

\ignore{
\input BoxedEPS
\SetRokickiEPSFSpecial
\HideDisplacementBoxes
\SetEPSFDirectory{./graphics/} }

\ignore{
\input BoxedEPS
\SetTexturesEPSFSpecial
\HideDisplacementBoxes
\SetEPSFDirectory{:graphics:}
}

\begin{document}

\title[Braid rigidity for path algebras]
{Braid rigidity for path algebras}

\author{Lilit Martirosyan }
\address{L.M. Department of Mathematics and Statistics\\ University of North Carolina\\ Wilmington \\North Carolina}
\email{martirosyanl@uncw.edu}

\author{Hans Wenzl}
\address{H.W. Department of Mathematics\\ University of California\\ San Diego,
California}

\email{hwenzl@ucsd.edu}

\begin{abstract}
Path algebras are a convenient way of describing decompositions of tensor powers
of an object in a tensor category. If the category is braided, one obtains representations
of the braid groups $B_n$ for all $n\in \N$.  We say that such representations are rigid if they are determined
by the path algebra and the representations of $B_2$. We show that besides the known
classical cases also the braid representations for the path algebra for the 7-dimensional
representation of $G_2$ satisfies the rigidity condition, provided $B_3$ generates 
$\End(V^{\otimes 3})$. We obtain a complete classification of ribbon tensor categories
with the fusion rules of $\g(G_2)$ if this condition is satisfied.
\end{abstract}
\maketitle

We say that a simple object $V$ in a semisimple tensor category $\Ca$ has the multiplicity 1
property if its tensor product with any simple object in $\Ca$ decomposes into a direct
sum of mutually non-isomorphic simple objects. This provides a convenient canonical
decomposition of $V^{\otimes n}$ into a direct sum of simple objects which are
labeled by paths, and a description of $\End(V^{\otimes n})$ via
path algebras.  E.g. for $V$ the vector representation of $Gl(N)$, the paths 
correspond to standard tableaux of certain Young diagrams.
If the category is braided, we also obtain representations of the braid group $B_n$
with respect to a basis labeled by these paths. 
For $Gl(N)$, we would obtain Young's orthogonal representations of
the symmetric groups. 

It is a fundamanental problem to classify all possible tensor categories for a given set
of tensor product (or fusion) rules. For braided tensor categories, an important tool
is to classify the corresponding representations of the braid group. This proved to be
successful in classifying braided tensor categories whose fusion rules were the ones
of a classical Lie group, see \cite{KW}, \cite{TW2}, where the braid representations
could be described in terms of Hecke algebras and $BMW$ algebras. 
Unfortunately, there does not seem to be a convenient algebraic description 
(via relations) of braid representations appearing for exceptional Lie groups.
This motivated our approach via path algebras. 

One can abstractly define braid representations compatible with path algebras,
see Eq \ref{braid2} and \ref{pathrep}.
We say that a path algebra is {\it braid rigid}, if any compatible braid representation
(with some mild additional conditions, see Definition \ref{braidrigiditydef})
is already uniquely determined by the image of the first braid generator.
This is the case for path algebras associated to classical Lie types. 
The main result in this
paper states that also the path algebra associated to the 7-dimensional
representation $V$ of $G_2$ is braid rigid, provided that the image of 
the braid group $B_3$
generates $\End(V^{\otimes 3})$. In particular, in these cases we obtain
the same braid representations as for the quantum group $U_q\g(G_2)$ 
of Lie type $G_2$.
This implies the classification of ribbon tensor categories whose fusion rules are the
ones of $\g(G_2)$ if they satisfy the condition about $\End(V^{\otimes 3})$ just stated. 
This result has already appeared before in \cite{MPS}, but our proof is
quite different and does not use any computer calculations.

Our approach was inspired by our previous work \cite{MW} where we
gave another proof that the braid groups generate $\End(V^{\otimes n})$
for $V$ the 7-dimensional representation $V$ of $U_q\g(G_2)$ for $q$
not a root of unity. We did this by finding quite explicit formulas for 
the path representation of braid groups for certain types of paths. Our main new
result in this paper is that the path representations for an abstract semisimple 
ribbon tensor category with the fusion rules of $\g(G_2)$ have to be isomorphic 
to the ones for the quantum group $U_q\g(G_2)$ for $q$ not a root of unity, 
at least when the condition  for the third tensor power of $V$
is satisfied.

Here is the content of our paper in more detail. We  review basic
definitions concerning path algebras and braided tensor categories in
the first section. We also define braid rigidity for path algebras there.
We then give the necessary combinatorial and algebraic information
about the Lie algebra $\g(G_2)$ in the second section. 
Let $V$ be the object in an abstract tensor category $\Ca$ of type $G_2$ corresponding
to the smallest nontrivial representation of $\g(G_2)$.
Then we show in the third section that if $B_3$ generates $\End_\Ca(V^{\otimes 3})$,
this braid representation has to be isomorphic to the corresponding
representation in $\End_\UU(V^{\otimes 3})$, where $\UU=Rep(U_q\g(G_2))$ for 
some $q$ not a root of unity.
The main technical result is then proved in the last section. We show that
the result in Section 3 can be extended to all tensor powers of $V$ by 
proving that the corresponding braid reprseentations are path rigid.
 We use this to
classify all ribbon tensor categories $\Ca$ with the fusion rules of $G_2$,
subject to the already mentioned condition concerning $\End_\Ca(V^{\otimes 3})$.

$Acknowledgment:$ Both authors would like to thank UNC Wilmington for its support through the Cahill Grant, and the Max Planck Institute in Bonn for its hospitality and support.

\section{Basic Definitions}

\subsection{Path algebras}\label{Litpath}
The notion of path algebras (or equivalent versions of it)  has been known
in many contexts such as operator algebras, representation theory
and algebra for a long time. We review some basic facts here, which will also
help to fix notations.

Let $\La$ be a set of labels
with distinguished label $0$
together with a  not necessarily symmetric relation $\rightarrow$. 
A path of length $n$ is a map $t:\{ 0,1,\ ...\ n\}\to\Lambda$ such that
$t(0)=0$ and $t(i)\rightarrow t(i+1)$ for $0\leq i<n$.
We denote by $\Pn$  the set of all paths of length $n$.
 We define algebras $C_n$ by
\begin{equation}\label{Cndef}
 C_n\ \cong\ \bigoplus_\nu M_{m(\nu, n)},
\end{equation}
where $m(\nu, n)$ is the number of paths $t$ of length $n$ with $t(n)=\nu$
and $M_m$ are the $m \times m$ matrices.   Let $W(\nu,n)$ be
a simple $C_n$-module labeled by the label $\nu$. It follows from the definitions
that it has a basis labeled by the paths in $\Pn$ which end in $\nu$.
Its decomposition into simple $C_{n-1}$ modules is given by the map $t\mapsto t'$,
where $t'$ is the restriction of $t$ to $\{ 0, 1,,\ ..., n-1\}$.
Hence we have the following isomorphism of $C_{n-1}$-modules:
\begin{equation}\label{restriction}
W(\nu,n) \cong \oplus_\mu  W(\mu, n-1),
\end{equation}
where $\mu$ runs through all labels $\mu$ which are endpoint of a path of length $n-1$
such that $\mu\rightarrow \nu$.

\begin{definition}\label{pathalgdef}
The path algebra $\P$ corresponding to the label set $\Lambda$ with 
the relation $\rightarrow$ is given by the sequence of algebras $C_n$
with the embeddings $C_{n-1}\subset C_n$ defined by \ref{restriction}.
\end{definition}

\begin{example}\label{pathexample}
The standard example for a path algebra is given by the labeling set
$\Lambda$ consisting of all Young diagrams with 0 being the empty Young
diagram, and $\mu\rightarrow \nu$ if $\mu\subset \nu$ and $|\nu|=|\mu|+1$,
i.e. $\nu$ is obtained by adding a box to $\mu$ . Then each path
corresponds to a Young tableau, and $C_n\cong \C S_n$.
We refer to this path algebra as the path algebra from Young's lattice.
\end{example}

\begin{remark} Our version is not the most general version of path algebra.
It will be clear from the next subsection that there also exist interesting
examples where one replaces our relation $\rightarrow$ by 
nonnegative integers $k(\mu,\nu)$ for each ordered pair $(\mu,\nu)\in \Lambda\times\Lambda$.
As this more complicated case will not be relevant in our paper, we stick
to this simpler version.
\end{remark}

\subsection{Tensor categories}
In this and the following subsection $\Ca$ will denote a semisimple
tensor category, whose $Hom$ spaces are complex vector spaces.
The reader not familiar with tensor categories can safely think of
$\Ca$ being the representation category of a Drinfeld-Jimbo quantum
group, or just of the corresponding semisimple Lie algebra.
Let $\Lambda$ be a labeling set for the simple objects of $\Ca$,
where $0$ is the label for the trivial object.
We also assume that $V$ is a simple object of $\Ca$ with the
{\it multiplicity 1 property}, i.e
$V_\mu\otimes V$ is a direct sum of mutually non-isomorphic simple
 objects, for any simple object $V_\mu$ in $\Ca$.
We then define the relation $\mu\rightarrow \nu$ if
$V_\nu\subset V_\mu\otimes V$.
This allows us to give a
fairly simple description of $\End_\Ca( V^{\otimes n})$ via paths.

\begin{theorem}\label{ldecomposition} We have a direct sum
decomposition of objects in $\Ca$ given by
$$ V^{\otimes n}=\bigoplus_\nu m(\nu,n)V_\nu,$$
where the multiplicity $m(\nu,n)$ is given by the number of paths
in $\Pn$ which end in $\nu$. In particular, we have
\begin{equation}\label{Cndef}
 C_n=\bigoplus_\nu M_{m(\nu, n)}\ \cong\ \End_\Ca( V^{\otimes n})
\end{equation}
where $M_k$ are the $k\times k$ matrices.
\end{theorem}

\begin{remark} As indicated at the end of the previous
subsection, it is not hard to give a path algebra description of
$\End_\Ca(V^{\otimes n})$ also if $V$ does not have the multiplicity
1 property. This can be done in terms of 
Littelmann paths (see \cite{Li}). The simpler version here has been
known much longer, see e.g. \cite{St} and references there.
\end{remark}

\begin{corollary}\label{pathc}
There exists an assignment $t\in \Pn\mapsto p_t\in C_n=
\End_\Ca( V^{\otimes n})$
such that $p_tV^{\otimes n}$ is an irreducible $\Ca$-object 
labeled by $t(n)$, and such that $p_tp_s=\delta_{ts}p_t$.
The idempotents $p_t$ are uniquely defined by the properties above
and the following one: If $s\in \P_{n-1}$, we have
$$p_s\otimes id=\sum_{t,\ t'=s} p_t.$$
\end{corollary}

\subsection{Path representations} We denote by $\Pn(\nu)$ 
all paths of length $n$ in $\Pn$ which end in $\nu$. 
One checks easily that $z_\nu^{(n)}=\sum_{t\in \Pn(\nu)}p_t$ is
a central idempotent in $C_n=\End_\Ca( V^{\otimes n})$.

By definition, we can define a basis $(v_t)_{t\in\Pn(\nu)}$
for the simple $C_n(\nu)$-module $W(\nu,n)$; often we will just write $ t$
for $v_t$.
Here the vector $v_t$ spans the image
of $p_t$ for each $t\in \Pn(\nu)$ and it is uniquely determined
up to scalar multiples.
Let $\delta, \nu$ be dominant weights for which
$V_\delta\subset  V^{\otimes n-k}$ and $V_\nu\subset  V^{\otimes n}$,
and let $\P_k(\delta,\nu)$ be the set of all paths of length $k$
from $\delta$ to $\nu$, with paths as defined
in Section \ref{Litpath}. Let $W_k(\delta,\nu)$ be the vector space
spanned by these paths and let $t$ be a fixed path in $\P_{n-k}(\delta)$.
Then we obtain a representation of
$\End_\Ca(V^{\otimes k})$ on $W_k(\delta,\nu)$ by
\begin{equation}\label{skew}
a\in \End_\Ca(V^{\otimes k}) \mapsto (p_t\otimes a)z_\nu^{(n)};
\end{equation}
here we used the obvious bijection between elements $s\in \P_k(\delta,\nu)$
and paths $\tilde s\in \P_n(\nu)$
for which $\tilde s_{|[0,n-k]}=t$, i.e. $\tilde s$ is the extension of $t$ by $s$.

\subsection{Braided tensor categories}\label{braidedtensor}  
We recall a few basic facts about braided and ribbon
tensor categories, see e.g. \cite{Ks}, \cite{Turaev} for more details.
This serves mostly as motivation for the definitions in the next subsection.
A braided tensor category $\Ca$ has canonical isomorphisms
$c_{V,W}: V\otimes W\to W\otimes V$ for any
objects $V$, $W$ in $\Ca$. They satisfy the condition
\begin{equation}\label{braid1}
c_{U, V\otimes W}\ =\ (1_V\otimes c_{U,W}) (c_{U,V}\otimes 1_W),
\end{equation}
and a similar identity for $c_{U\otimes V, W}$. Let $B_n$ be Artin's braid groups, 
given by generators $\sigma_i,
\ 1\leq i\leq n-1$ and relations $\sigma_i\sigma_{i+1}\sigma_i=
\sigma_{i+1}\sigma_i\sigma_{i+1}$ as well as $\sigma_i\sigma_j=\sigma_j
\sigma_i$ for $|i-j|\geq 2$.
One can show that we obtain a representation of the braid group $B_n$
into $\End(V^{\otimes n})$ for any object $V$ in $\Ca$ via the map
\begin{equation}\label{braid2}
\sigma_i\ \mapsto 1_{i-1}\otimes  c_{V,V} \otimes 1_{n-1-i},
\end{equation}
where $1_k$ is the identity morphism on $V^{\otimes k}$.
Using a path basis $(t)$ as in the last subsection,
we can express the action of $\sigma_i$ via a matrix $A_i$ such that
\begin{equation}\label{pathrep}
\sigma_i\mapsto A_i: t\to \sum_s a_{st}^{(i)} s,
\end{equation}
where the summation goes over paths $s$ for which $s(j)=t(j)$ for $j\neq i$;
this follows from Eq \ref{skew} with $n=i+1$ and $k=2$.
As the vectors $t$s are uniquely determined up to rescaling,
it also follows that the matrix entries of $A_i$ are uniquely determined up
to conjugation by a diagonal matrix.

An associated ribbon braid structure is given by maps $\Theta_W: W\to W$ satisfying
\begin{equation}\label{ribbon}
\Theta_{V\otimes W}\ =\ c_{W,V}c_{V,W} (\Theta_V\otimes \Theta_W).
\end{equation}
Let $\Delta_n\in B_n$ be defined inductively by $\Delta_2=\sigma_1$ and
$\Delta_n=\Delta_{n-1}\sigma_{n-1}\sigma_{n-2}\ ...\ \sigma_1$. 
Then it is well-known that $\Delta_n^2=(\sigma_1\sigma_2\ ...\ \sigma_{n-1})^n$
generates the center of $B_n$. One can then prove by induction on $n$
that 
\begin{equation}\label{twistbraid}
\Theta_{V^{\otimes n}}= \Delta_n^2 \Theta_V^{\otimes n}.
\end{equation}
 If $V_\la$ is a simple object, 
the ribbon map just acts via a scalar, which we will denote by $\Theta_\la$.  
Let us also assume that the representation of $B_n$ into $\End(V^{\otimes n})$
is semisimple. Then the central element $\Delta_n^2$ acts in the simple
component labeled by $\alpha$ via a scalar denoted by $z_{\al, n}$.
If the $B_n$-representation labeled by $\alpha$ acts nontrivially on the
$\End(V^{\otimes n})$-module $W_\la^{(n)}$, then it follows
from Eq \ref{twistbraid} that
\begin{equation}\label{scalarmatch}
\Theta_\la\ =\ z_{\al, n}\Theta_v^n,
\end{equation}
where we identified $\Theta_V$ with the scalar via which it acts on $V$.

\subsection{Braid rigidity}\label{brigid} We now translate the notions
of the previous section into the language of path algebras.

\begin{definition}\label{braidrigiditydef} Let $\P$ be a path algebra.

(a) We call a system of representations of braid groups 
{\it representations of type $\P$} if
the braid generators act on paths as in \ref{pathrep}. Moreover, we also 
require that the central element $\Delta_n^2\in B_n$ acts via a fixed scalar
$z_{\la,n}$ on every path of length $n$ which ends in $\la$.

(b) We call a path algebra $\P$ {\it braid rigid} if any non-trivial
braid representation of
type $\P$ is uniquely determined by the image of $\sigma_1$; see the example
below for the description of trivial braid representations of type $\P$.
\end{definition}

\begin{example} 1. We can always define trivial braid representations
for any path algebra $\P$ as follows. We fix a non-zero number $\al$
and we assign to each path $s$ of length 2
an eigenvalue $\al_s=\pm \al$ of $\sigma_1$. Then we define the action of
$\sigma_i$ on a path $t$ to be equal to multiplication by $\al_s$ if the
restriction of $t$ to $\{ 0, 1, 2\}$ is equal to $s$. It is easy to check that
we obtain a path representation of $B_n$ for any $n$ which is a direct
sum of 1-dimensional abelian representations.

2. Let $\P$ be the path algebra given by Young's lattice.
We claim that the coresponding path algebra is braid rigid provided that the
ratio between the two eigenvalues of $\sigma_1$ is not a root of  unity.
Indeed, as we only have two paths of length two, the image of
$\sigma_1$ has at most two eigenvalues.
It is well-known that in this case we obtain representations of the Hecke algebras 
$H_n(q)$ of type $A_{n-1}$,
where $q=-\al_1/\al_2$ for $\al_1$ and $\al_2$ being the eigenvalues
of the image of $\sigma_1$. If the braid representation is nontrivial, 
the representation of $B_3$ on the 2 paths ending in the Young diagram $[21]$
has to be irreducible. Using the $q$-Jucys-Murphy approach,
see e.g. \cite{MW}, Lemma 1.8 for a review, we can then inductively
compute all matrix entries (up to rescaling of basis vectors)
for any path representation.

3. The same statement is also true if $-\al_1/\al_2$ is a primitive
$\ell$-th root of unity if we restrict the label set $\Lambda$ to the set
$\Lambda^{(k,\ell)}$ of so-called $(k,\ell)$-diagrams, i.e. to Young diagrams
with $\leq k$ rows such that $\la_1-\la_k\leq \ell-k$, see \cite{WHecke}
for details.

4. One can similarly also show that the path representations for the path algebra
generated by the vector representation $V$ of an orthogonal or symplectic
group are braid rigid. This follows essentially from \cite{TW2}, where a complete
classification of braided tensor categories was given for which the fusion rules
are the ones of the representation category of an orthogonal or symplectic group.
The main point of the proof there was to show that $\End_\Ca(V^{\otimes n})$ was
given by a quotient of the so-called $BMW$-algebra, see \cite{BBMW}, \cite{Mu}.
\end{example}

\subsection{Matrix blocks} We again assume $\Ca$ to be a general ribbon
tensor category as in Section \ref{braidedtensor}.
It follows from Eq \ref{skew} and \ref{pathrep} 
that the matrix $A_{n-1}$
acts in blocks leaving invariant path spaces $W_2(\delta, \la)$
spanned by a basis $(v_t)$ labeled by all
paths of length 2 from $\delta$ to $\lambda$.
It follows from the definitions that 
\begin{equation}\label{Wdef}
 W_2(\delta, \la)\cong \Hom_\Ca(V_\lambda, V_\delta\otimes V^{\otimes 2}).
\end{equation}
Indeed, as the image of $\sigma_{n-1}$ commutes with 
$\End(V^{\otimes n-2})$, the only relevant part of the path $t$
for the action of $A_{n-1}$ are the weights $t(i)$, $n-2\leq i\leq n$.
We will consider certain cases in which we can calculate the matrix
entries of $A_{n-1}$. More precisely we
consider the following cases:
\renewcommand{\labelenumi}{\alph{enumi})}
\begin{enumerate}
\item We have $\dim\ W_2(\delta,\la)=2$, and $A_{n-1}$ acts with two distinct eigenvalues on it.
\item The action of $A_{n-1}$ is diagonalizable on $W_2(\delta,\la)$ with exactly three distinct eigenvalues,
with one of them having multiplicity 1. We define $q$ such that the ratio of the other two eigenvalues
is equal to $-q^2$.
\end{enumerate}
The following proposition is a reformulation of results in \cite{Wexc}
and \cite{MW}:

\begin{lemma}\label{AB2representations} Assume that $A_{n-1}$ and $W_2(\delta,\la)$
satisfy the conditions just stated. Then its entries with respect to the path basis of
$W_2(\delta,\la)$ can be calculated in terms of the eigenvalues of $A_{n-1}$ and the 
entries of the rank 1  eigenprojection $P$, up to conjugation by a diagonal matrix. In particular,
for each such block all the off-diagonal matrix entries of $A_{n-1}$ are nonzero
if the corresponding entries of $P$ are nonzero.
\end{lemma}

$Proof.$ This is a consequence of \cite{MW}, Lemma 1.8,  Prop. 1.6, Lemma 1.8
and Lemma 3.2. We give some details for the reader's convenience.
 If $\dim W_2(\delta,\lambda) = 2$, the claim follows from a well-known
$q$-version of the Jucys-Murphy approach, see e.g. \cite{MW},  Lemma 1.8
for details.  If $A_{n-1}$ has three eigenvalues, let $P$ be the eigenprojection 
of $A_{n-1}$ for the eigenvalue with multiplicity 1.
Moreover, let $q^\al$ be a scalar such that $A'=q^{e_\al}A_{n-1}$ has
eigenvalues $q$, $-q^{-1}$ and $r^{-1}$ such that $A'P=r^{-1}P$. 
As, by construction we have
$$A'-(A')^{-1}\ =\ (q-q^{-1})I-(r-r^{-1}+q-q^{-1})P,$$
we can calculate the matrix entries of $A'$ from the equation
$$(1-q^{e(t)+e(s)})a'_{ts}=(q-q^{-1})\delta_{ts}-(r-r^{-1}+q-q^{-1})p_{ts},$$
where $q^{e(t)}$ is the scalar via which $\Delta_n^2\Delta_{n-1}^{-2}$ acts
on the path $t$,
see \cite{Wexc}, Lemma 4.1 for details.

\ignore{
\renewcommand{\labelenumi}{\alph{enumi})}
\begin{enumerate}
\item They satisfy the braid relation $ATAT=TATA$,
\item The matrix $A$ satisfies the relation $A-A^{-1}=(q-q^{-1})(1-mP)$, where $m=1 + (r-r^{-1} )/(q-q^{-1})$, 
and where $P$ is a rank 1 eigenprojection of $A$.
\item The central element $TATA$ acts as the identity on $W$.
\item We assume that $T$ is a diagonal matrix with eigenvalues $q^{e(t)}$
    where $t$ runs through a labeling set for a basis of $W$.
\end{enumerate}
}

\vskip .3cm
\ignore{
\begin{proposition}\label{AB2representations}
(a) The matrix entries of $A$ and $P$ are related by the equation
$$(1-q^{e(t)+e(s)})a_{ts}=(q-q^{-1})\delta_{ts}-(r-r^{-1}+q-q^{-1})p_{ts}.$$
(b) The diagonal entry $d_s=p_{ss}$ is equal to zero only if $e(s)=\pm 1$.
(c) The matrix $A$ is uniquely determined up to conjugation by an invertible
diagonal matrix.
\end{proposition}}

\section{The example $G_2$}

\subsection{Quantum groups}\label{quantum} The best known examples
of ribbon categories are given by the representation categories $\UU=Rep(\U)$
of a Drinfeld-Jimbo quantum group $\U=U_q\g$, where $\g$ is a semisimple
Lie algebra.
We assume as ground ring the field
$\C(q)$ of rational functions in the variable $q$. It is well-known
that in our setting the category $\UU$
of integrable representations of $\U$ is semisimple,
and it has the same Grothendieck semiring
as the original Lie algebra. We shall need the following result
due to Drinfeld \cite{Dr}.
\begin{proposition}
\label{prop:eigenvalue}
Let $V_\la, V_\mu, V_\La=V$ be
simple $\U$-modules with highest weights
$\la, \mu, \La$ respectively, and such that $V_\mu$ is a submodule of
$V_\la\otimes V_\La$. Let us write $c_{\la,\mu}$ for the braiding morphism
$V_\la\otimes V_\mu\to V_\mu\otimes V_\la$. Then
$$(c_{{\la},{\La}}c_{{\La},{\la}})_{|V_\mu}=
q^{C_\mu- C_\la-C_\La}1_{V_\mu},$$
where for any weight $\gamma$ the quantity $C_\gamma$ is given by
$C_\gamma=(\gamma +2\rho ,\gamma)$. Here, $\rho$ is the Weyl vector.
Moreover, the twisting factors $\Theta_\la$ are given by
$$\Theta_\la\ =\ q^{C_\la}.$$
\end{proposition}
\noindent

\subsection{Path algebra for $G_2$}\label{G2basics} We will be particularly interested in the case
with $\g=\g(G_2)$ and $V$ its simple 7-dimensional representation.
We first recall some basic facts about its roots and weights (see e.g. \cite{humphreys}, \cite{Kc}).

With respect to the orthonormal unit vectors $\eps_1$, $\eps_2$, $\eps_3$ of $\R^3$, 
the roots of $\g$ can be written $\Phi=\pm \{\eps_1-\eps_2,\eps_2-\eps_3, \eps_1-\eps_3, 2\eps_1-\eps_2-\eps_3, 2\eps_2-\eps_1-\eps_3, 2\eps_3-\eps_1-\eps_2\}$. 
The base can be chosen $\Pi=\{\al_1=\eps_1-\eps_2, \al_2=-\eps_1+2\eps_2-\eps_3\}$. The Weyl vector is given by $\rho=2\eps_1+\eps_2-3\eps_3$ 
and the Weyl group is $D_6$. The fundamental dominant weights are given by $\{\La_1=\eps_1-\eps_3, \La_2=\eps_1+\eps_2-2\eps_3 \}$.
The following describes the dominant Weyl chamber:

\begin{center}
  \begin{tikzpicture}[scale=.4]
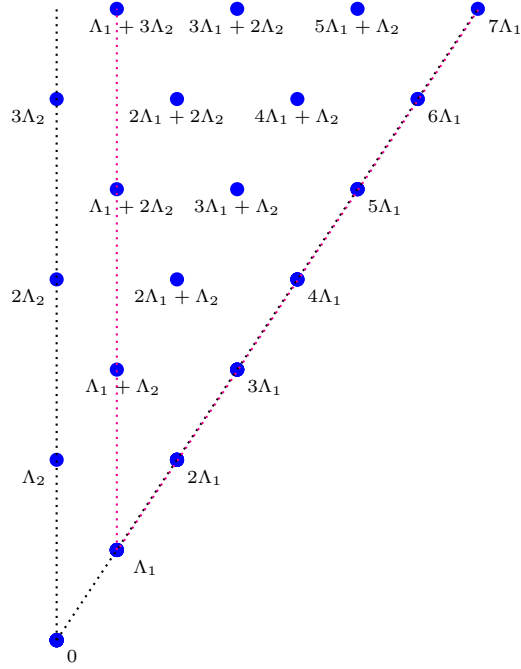
    \draw (0,-3) node[anchor=east]  {Weyl Chamber for $\g$};
   \foreach \y in {0,...,7} \foreach \x in {0,...,\y}
    \draw[xshift=\x cm,thick,blue,fill=blue](\x cm,3*\x cm) circle (.2cm);
       \foreach \x in {0,...,5}
    \draw[xshift=\x cm,thick,blue,fill=blue] (\x cm,3*\x cm +6 cm) circle (.2cm);
        \foreach \x in {0,...,3}
    \draw[xshift=\x cm,thick,blue,fill=blue] (\x cm,3*\x cm +12 cm) circle (.2cm);
     \foreach \x in {0,...,1}
    \draw[xshift=\x cm,thick,blue,fill=blue] (\x cm,3*\x cm +18 cm) circle (.2cm);
    \foreach \x in {2,...,7}
    \draw[xshift=\x cm] (\x cm,3*\x cm) node [anchor=north west]{\tiny $\x\La_1$};
      \foreach \x in {2,...,5}
    \draw[xshift=\x cm] (\x cm,3*\x cm + 6 cm) node [anchor=north]{\tiny $\x\La_1+\La_2$};
     \foreach \x in {2,...,3}
    \draw[xshift=\x cm] (\x cm,3*\x cm + 12 cm) node [anchor=north]{\tiny $\x\La_1+2\La_2$};
        \foreach \y in {2,...,3}
    \draw[yshift=\y cm] (0cm, 5*\y cm) node [anchor=north east]{\tiny $\y\La_2$};
         \foreach \y in {2,...,3}
    \draw[yshift=\y cm] (2.5 cm, 5*\y cm+ 3cm) node [anchor=north]{\tiny $\La_1+\y\La_2$};
   	
	  \draw (0 cm, 0 cm) node[anchor=north west]  {\tiny $0$};
       \draw (0 cm, 6 cm) node[anchor=north east]  {\tiny $\La_2$};
       \draw (2.2 cm, 3 cm) node[anchor=north west]  {\tiny $\La_1$};
       \draw (2.2 cm, 9 cm) node[anchor=north]  {\tiny $\La_1+\La_2$};
 	  \draw[dotted,thick] (0,0) -- +(14, 21);
 	  \draw[dotted,thick] (0,0) -- +(0, 21);
 	 	  \draw[dotted,thick, magenta] (2,3) -- +(0, 18);
 	  \draw[dotted,thick, magenta] (2.1,3.1) -- +(12, 18);


  \end{tikzpicture}
\end{center}
It follows from Weyl's dimension formula that the $q$-dimension of the $U_q\g(G_2)$-module $V_\la$ 
with highest weight  $\la=(\la_1,\la_2, -\la_1-\la_2)$ 
is equal to
\begin{equation}\label{dimensionformula}
\dim_q V_\la\ =\ \frac{[\la_1-\la_2+1][2\la_1+\la_2+5][\la_1+2\la_2+4][3\la_1+6][3\la_2+3][3(\la_1+\la_2)+9]}
{[1][5][4][6][3][9]}
\end{equation}

\subsection{Tensor product rules}\label{tensorsect}
We will study the path representations with respect to the smallest nontrivial
$\g$-module $V=V_{\Lambda_1}$. We will need to know how to tensor irreducible
representations with $V$.
We will review this here for the reader's convenience.
The representation $V$ has dimension 7, with its weights being
the short roots of $\g$ together with the zero weight.
The decomposition of the tensor product 
\begin{equation}\label{tensorrule}
V_\la\otimes V\ \cong\ \bigoplus_\mu V_\mu 
\end{equation}
with $V_\la$ a simple
module with highest weight $\la=a\Lambda_1+b\Lambda_2$ can
be described as follows (see e.g. \cite{MW}, Prop. 2.1 and Remark 2.2): 
Consider the hexagon centered at $\la$
and with corners $\la+\om$, with $\om$ running through
the short roots of $\g$. If this hexagon is contained in the dominant
Weyl chamber $C$, then $V_\la\otimes V$ decomposes into the direct sum
of irreducibles $\g$-modules whose highest weights are given
by the corners and the center of
the hexagon. If it is not contained in $C$, leave out all
the corners of the hexagon which are not in $C$; moreover, if
$\la=b\La_2$, also leave out $\la$ itself.
Using this, we can draw the Bratteli diagram for $V^{\otimes n}$.

  \begin{tikzpicture}[scale=.4]

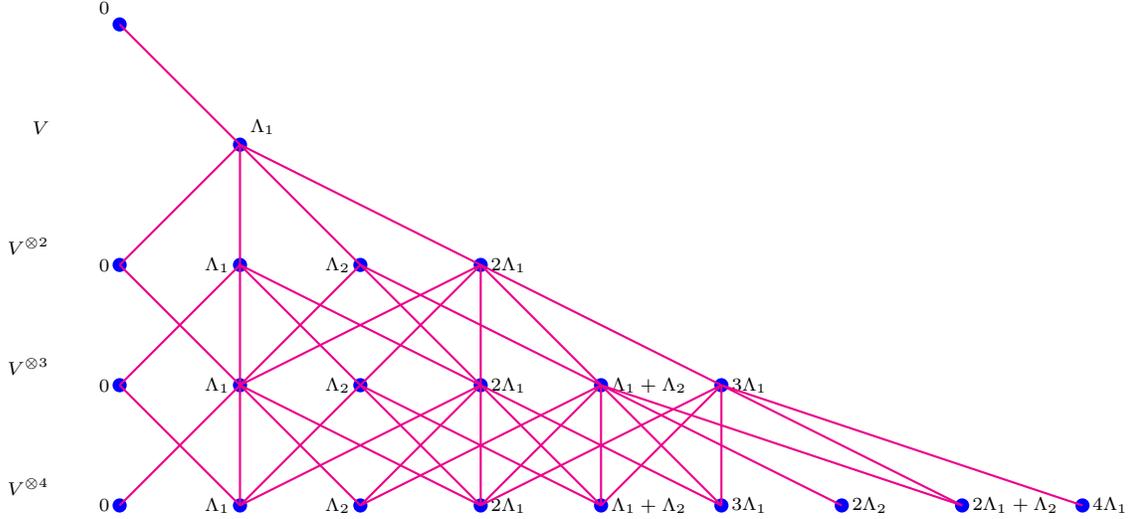
    \draw (0,-16) node[anchor=east]  {Bratteli diagram for $V^{\otimes n}$};
        \draw (-10, 4) node[anchor=south east]  {\tiny $0$};

    \draw [thick,blue,fill=blue](-10 cm,4 cm) circle (.2cm);

     \draw (-6, 0) node[anchor=south west]  {\tiny $\La_1$};
     \draw (-12, 0) node[anchor=south east]  {\tiny $V$};
\foreach \x in {2,...,4}      \draw (-12, -4*\x+4) node[anchor=south east]  {\tiny $V^{\otimes \x}$};

\foreach \x in {1,...,3}     \draw (-10, -4*\x) node[anchor= east]  {\tiny $0$};
\foreach \x in {1,...,3}     \draw (-6, -4*\x) node[anchor= east]  {\tiny $\La_1$};
\foreach \x in {1,...,3}     \draw (-2, -4*\x) node[anchor= east]  {\tiny $\La_2$};
\foreach \x in {1,...,3}     \draw (2, -4*\x) node[anchor= west]  {\tiny $2\La_1$};
\foreach \x in {2,3}     \draw (6, -4*\x) node[anchor= west]  {\tiny $\La_1+\La_2$};
\foreach \x in {2,3}     \draw (10, -4*\x) node[ anchor= west]  {\tiny $3\La_1$};
\draw (14, -12) node[anchor= west]  {\tiny $2\La_2$};
 \draw (18, -12) node[anchor= west]  {\tiny $2\La_1+\La_2$};
  \draw (22, -12) node[anchor= west]  {\tiny $4\La_1$};

    \draw [thick,blue,fill=blue](-6 cm,0 cm) circle (.2cm);
 \foreach \x in {0,...,3}
   \draw[xshift=\x cm,thick,blue,fill=blue] (3*\x cm-10 cm,-4 cm) circle (.2cm);

 \foreach \x in {0,...,5}
   \draw[xshift=\x cm,thick,blue,fill=blue] (3*\x cm-10 cm,-8 cm) circle (.2cm);
\foreach \x in {0,...,8}
\draw[xshift=\x cm,thick,blue,fill=blue, outer sep=0.3pt,
    inner sep=0.5pt] (3*\x cm-10 cm,-12 cm) circle (.2cm);

      \draw[thick, magenta] (-10,4) -- +(4, -4);

    	  \draw[thick, magenta] (-6,0) -- +(8, -8);
    	  \draw[thick, magenta] (-6,0) -- +(16, -8);
    	  \draw[thick, magenta] (-6,0) -- +(0, -12);
	  \draw[thick, magenta] (-10,-4) -- +(4, 4);
	  \draw[thick, magenta] (-6,-4) -- +(8, -4);
	  \draw[thick, magenta] (-6,-4) -- +(4, -4);

	  \draw[thick, magenta] (-6,-8) -- +(8, -4);

    	  \draw[thick, magenta] (-10,-4) -- +(8, -8);
    	  \draw[thick, magenta] (-10,-8) -- +(4, -4);
    	  \draw[thick, magenta] (-10,-8) -- +(4, 4);
    	  \draw[thick, magenta] (-10,-12) -- +(8, 8);
    	  \draw[thick, magenta] (-2,-8) -- +(4, -4);
    	  \draw[thick, magenta] (-2,-8) -- +(8, -4);
    	  \draw[thick, magenta] (-2,-4) -- +(8, -4);
    	  \draw[thick, magenta] (-2,-8) -- +(-4, -4);
    	  \draw[thick, magenta] (2,-4) -- +(0, -8);
    	  \draw[thick, magenta] (2,-4) -- +(-8, -4);
    	  \draw[thick, magenta] (2,-4) -- +(4, -4);
    	  \draw[thick, magenta] (2,-8) -- +(-8, -4);
    	  \draw[thick, magenta] (2,-8) -- +(8, -4);
    	  \draw[thick, magenta] (2,-8) -- +(4, -4);
    	  \draw[thick, magenta] (2,-8) -- +(-4, -4);
    	  \draw[thick, magenta] (2,-4) -- +(-4, -4);
    	  \draw[thick, magenta] (6,-8) -- +(-4, -4);
    	  \draw[thick, magenta] (6,-8) -- +(0, -4);
    	  \draw[thick, magenta] (6,-8) -- +(-8, -4);
    	  \draw[thick, magenta] (6,-8) -- +(4, -4);
    	  \draw[thick, magenta] (6,-8) -- +(8, -4);
    	  \draw[thick, magenta] (6,-8) -- +(12, -4);
    	  \draw[thick, magenta] (10,-8) -- +(0, -4);
    	  \draw[thick, magenta] (10,-8) -- +(12, -4);
    	  \draw[thick, magenta] (10,-8) -- +(-4, -4);
    	  \draw[thick, magenta] (10,-8) -- +(8, -4);
    	  \draw[thick, magenta] (10,-8) -- +(-8, -4);

\end{tikzpicture}

\subsection{Braid representations for $U_q\g(G_2)$}
We normalize the invariant product on the weight lattice such that
$(\La_1,\La_1)=2$ and $(\La_2,\La_2)=6$. With these conventions we get
the values $C_\nu=0, 28, 12, 24$ for $\nu=0,2\La_1,\La_1,\La_2$.
Hence  it follows from Proposition
\ref{prop:eigenvalue} that the eigenvalues of $R_{V,V}$ are given by $q^{-12}, q^{2},-q^{-6}$
and $-1$. So if $\alpha=q$, $A'=\alpha^{-1}R_{V,V}$ has the desired eigenvalues $q$ and $-q^{-1}$ 
for the representations $V_{2\La_1}$ and $V_{\La_2}$.
 As we shall see in a moment, it will be convenient to
associate with $P$ the eigenprojection of $A$ projecting onto
$V=V_{\La_1}\subset V^{\otimes 2}$, which corresponds to the eigenvalue $-q^{-6}$. Indeed, let
$W=\Hom_\U(V_\la, V_\delta\otimes V^{\otimes 2})$. Then $P$ is the projection onto
the subspace $\Hom_\U(V_\la, V_\delta\otimes V_{\La_1})$ of $W$,
given by the embedding $V\subset V^{\otimes 2}$.
As all weights of $V$ have multiplicity 1, the multiplicity of $V_\la$ in $V_\delta\otimes V_{\La_1}$
is at most 1.  Hence $P$ has at most rank 1 in $\Hom_\U(V_\la, V_\delta\otimes V^{\otimes 2})$.
In particular, the conditions for Lemma \ref{AB2representations} are satisfied.
We can now refine the results of that lemma in our setting as follows, restating results
which have already appeared before in \cite{Wexc} and \cite{MW}. Let 
$$t: \delta \rightarrow \mu_t \rightarrow \lambda$$
be a path of length 2 and let
$$e(t)=C_{\mu_t}-\frac{1}{2}(C_\lambda +C_\delta)+1,$$
where $C_\gamma = (\gamma +2\rho,\gamma)$ for a weight $\gamma$.

\begin{lemma}\label{smallblockss}
Consider the space $W_2(\delta,\lambda)$ with a basis labeled by paths of length 2
from $\delta$ to $\lambda$. If $\la\neq \delta$, then the matrix $A_{n-1}$ can be calculated
up to conjugation by a diagonal matrix  as in 
Lemma \ref{AB2representations}. In particular, all off-diagonal entries of $A_{n-1}$ are
well-defined and not equal to 0 if $q$ is not a root of unity.
\end{lemma}

$Proof.$ This result was essentially already shown in \cite{MW}.
We give a proof here for the reader's convenience. Let $P_\gamma$ be the eigenprojection of $A_{n-1}$
corresponding to the representation $V_\gamma\subset V^{\otimes 2}$.
It follows from the definitions that it acts on $_2W(\delta,\gamma)$ as a
rank $c_{\delta,\gamma}^{\lambda}$ idempotent, where 
$c_{\delta,\gamma}^{\lambda}$ is the multiplicity of $V_\lambda$ in 
$V_\delta\otimes V_\delta$. This rank is equal to 0 for $\gamma=0$
unless $\delta=\lambda$.  Hence $A_{n-1}$ can only act
with at most three distinct eigenvalues on $W_2(\delta,\gamma)$ for $\delta\neq\lambda$. 

It was shown in \cite{MW} that $A_{n-1}$ acts with two eigenvalues only if
the dimension of $W_2(\delta,\lambda)$ is equal to 2. It is well-known how to calculate
the matrix coefficients in this case via the Jucys-Murphy approach, see e.g. 
\cite{MW}, Lemma 1.8 for a review and a precise statement.

If $A_{n-1}$ acts with three distinct eigenvalues, the eigenprojection $P$ for $\gamma=\Lambda_1$
has rank 1. 
It was shown in \cite{MW}, Lemma 3.2 that its diagonal entry $d_s$ for the path $s$ is given by
$$d_s=\frac{[e(s)+1]}{1- [7]}\ \prod_{t\neq s}
\frac{[(e(s)+e(t))/2]}{ [(e(s)-e(t))/2]}.$$
It was shown in \cite{MW} Proposition 1.6 and Lemma 2.6 that these entries are nonzero
and well-defined for $q$ not a root of unity. As $P$ is a rank 1 idempotent, this also shows
that all of its matrix entries are nonzero. We can now calculate the matrix entries of $A_{n-1}$
as shown in the proof of Lemma \ref{AB2representations}. In particular, this shows that also
the off-diagonal entries of $A_{n-1}$ are well-defined and nonzero for $q$ not a root of unity.

\medskip
 The following result has first been shown in 
\cite{LZ}, with different proofs also given in \cite{Ms} and \cite{MW}:

\begin{theorem}\label{firstfund} (First Fundamental Theorem)
Let $V$ be the 7-dimensional representation of \U=$U_q(\g(G_2))$ with highest weight $\Lambda_1$.
Then $\End_\U(V^{\otimes n})$ is generated by the image of the braid group $B_n$
in  $\End_\U(V^{\otimes n})$ for $q$ not a root of unity.
\end{theorem}

\section{Tensor categories of type $G_2$}

In the rest of this paper, we let $\Ca$ be a semisimple rigid ribbon tensor category
of type $G_2$.
 By this we mean that its simple objects $X_\la$ are labeled by the dominant integral weights $\la$
of $G_2$, 
and the decomposition of tensor products of simple objects is given by the tensor product rules for $G_2$.
See e.g. \cite{Ks}, \cite{Turaev} for precise definitions of the other terms. We shall first study
the braid representations corresponding to small tensor powers of the object corresponding
to the 7-dimensional irreducible representation $V$. The main result of this section
is  that the eigenvalues
of the braid generators are forced to be the same as in the quantum group case whenever the
braid representations generate $\End_\Ca(V^{\otimes n})$ for $n=2,3$. 

\subsection{Preparations}

 We shall use properties of ribbon categories,
in particular Eq \ref{scalarmatch} to
find constraints for the eigenvalues of $c_{V,V}$. The following elementary lemma
will be useful:

\begin{lemma}\label{braideigenvalues} Let $W$ be a representation of $B_3$ of dimension
$m$ on which $\sigma_1$ acts with eigenvalues $\alpha_i$, $1\leq i\leq m$ and on which
$\Delta_3^2$ acts via the scalar $z_{3,m}$. Then we have
$$z_{3,2}=-(\al_1\al_2)^3,\quad z_{3,3}=(\al_1\al_2\al_3)^2, \quad z_{3,4}=\sqrt{\al_1\al_2\al_3\al_4}^3,$$
where there exist representations for both choices of the square root for $m=4$.
\end{lemma}

$Proof.$ Let us first assume we have a representation of $B_n$ acting on an $m$-dimensional
vector space such that $\Delta_n^2$ acts via the scalar $z_{n,m}$.
Calculating the determinant of the matrix representing 
$\Delta_n^2=(\sigma_1\sigma_2\ ...\ \sigma_{n-1})^{n(n-1)}$
in two different ways, we obtain
\begin{equation}\label{detscalar}
z_{n,m}^m=\det(\sigma_1)^{n(n-1)};
\end{equation}
This equation does not determine which $m$-th root of
the determinant we have to take for $z_{n,m}$.
But it was shown in \cite{TW} that for dimension $m\leq 5$  representations of
$B_3$ as in the statement are essentially obtained by their eigenvalues.
Using the explicit braid representations in \cite{TW}, one can check the claim
by a direct calculation.

\subsection{Calculations of eigenvalues}
We use the notations as in Section \ref{G2basics}.  The dominant integral weights $\la$ are
of the form $\la=(\la_1,\la_2,-\la_1-\la_2)$ with $\la_1\geq \la_2\geq 0$.
In the following we will just write $\la=(\la_1,\la_2)$ for brevity.
So the highest weight of the 7-dimensional simple representation $V$ of $G_2$
 is given by $\Lambda_1=(1,0)$. The second fundamental weight
is $\Lambda_2=(1,1)$. 
It follows from the decomposition of $V^{\otimes 2}$,
see Eq \ref{tensorrule}, or \cite{MW}, Example 2.3  that
$c_{V,V}$ has four eigenvalues corresponding
to the subrepresentations $\1=V_{(0,0)}$, $V=V_{(1,0)}$, $V_{(1,1)}$ and
$V_{(2,0)}$ respectively. We will refer to the eigenvalue belonging
to $V_\la\subset V^{\otimes 2}$ by $\al_\la$.
It follows from \ref{scalarmatch} for $n=2$ that 
\begin{equation}\label{second}
\Theta_\la\ =\ \al_\la^2\Theta_{(1,0)}^2.
\end{equation}
As $\Theta_{(0,0)}=1$ (which can be deduced from \ref{ribbon} and the braiding
axioms),
we obtain $1=\al_{(0,0)}^2\Theta_{(1,0)}^2$ and $\Theta_{(1,0)}=\al_{(1,0)}^2\Theta_{(1,0)}^2$ from
\ref{second}.  Hence we have
\begin{equation}\label{thetalow}
\pm \al_{(0,0)}\ =\ 1/\Theta_{(1,0)}\ =\ \al_{(1,0)}^2.
\end{equation}
Observe that the representations labeled by  $(2,0)$ and $(1,1)$ appear with multiplicity 
3 and 2  in $V^{\otimes 3}$. It follows from Lemma \ref{braideigenvalues} and Eq. \ref{scalarmatch}
that
$$\al_{(2,0)}^2\Theta_{(1,0)}^2\ =\ \Theta_{(2,0)}\ =\ (\al_{(1,0)}\al_{(1,1)}\al_{(2,0)})^2\Theta_{(1,0)}^3,$$
from which we deduce, together with \ref{thetalow}
\begin{equation}\label{al2}
\al_{(1,1)}^2=1.
\end{equation}
 Again using Lemma \ref{braideigenvalues} and Eq. \ref{second} for $\la=(1,1)$
we obtain
$$\al_{(1,1)}^2\Theta_{(1,0)}^2\ =\ \Theta_{(1,1)}\ =\ -(\al_{(1,0)}\al_{(2,0)})^3\Theta_{(1,0)}^3.$$
We deduce from this, again using $\al_{(1,0)}^2\Theta_{(1,0)}=1$,  that
\begin{equation}\label{al3}
\al_{(1,0)}\al_{(2,0)}^3=-1.
\end{equation}
We have almost proved the following proposition:

\begin{proposition}\label{eigenvalues} Assume that the image of 
 $B_n$ in $\End_\Ca(V^{\otimes n})$ generates these algebras for $n=2,3$.
Then these representations are isomorphic to the corresponding ones
appearing in $\UU=Rep(U_q(\g(G_2)))$ for some $q$.
\end{proposition}

$Proof.$ If we set $\al_{(2,0)}=q^2$, 
it follows $\al_{(1,0)}=-q^{-6}$ from \ref{al3}, $\al_{(1,1)}=\pm 1$ from  \ref{al2} 
and $\al_{(0,0)}=\pm q^{-12}$ from \ref{thetalow}. Using Lemma \ref{braideigenvalues}
for the representation on $W((1,0),3)=\Hom(V_{(1,0)},V^{\otimes 3})$,
for which the dimension $m_\la=4$, we obtain for the scalar $z_{(1,0),3}$
by which $\Delta_3^2$ acts that
$$(\al_{(0,0)}\al_{(1,0)}\al_{(1,1)}\al_{(2,0)})^{6/4}\ =\ z_{(1,0), 3}\ =\ \Theta_{(1,0)}^{-2}\ =\ \al_{(1,0)}^4 = q^{-24},$$
where we used Eq. \ref{scalarmatch}.
It follows that the product of the eigenvalues must be a power of $q$, i.e. the number
of minus signs for the eigenvalues must be even. Otherwise we would not get a power of $q$ from the radical.
This implies $\al_{(1,1)}\al_{(0,0)}=-q^{-12}$, which  forces  the eigenvalues to be
as in the statment, or as in the statement with opposite signs. Also observe that
this also shows that the square root of the determinant of $\sigma_1$ is the same
as in the quantum group case.

It only remains to show that the second option with $\al_{(1,1)}=1$ and $\al_{(0,0)}=-q^{-12}$
can not occur for a ribbon tensor category. We obtain braid representations
with such eigenvalues for the negative $R$-matrix in the quantum group case after substituting $q^2$ by $-q^2$.
As the trivial representation $\1$ appears in $V^{\otimes 3}$, where $V=V_{(1,0)}$, the
 negative $R$-matrix violates the braiding axioms as follows: Recall that the
$R$-matrix $R_V$ for $V^{\otimes 3}\otimes V^{\otimes 3}$ is the image of the braid 
$\sigma_3\sigma_2\sigma_1\sigma_4\sigma_3\sigma_2\sigma_5\sigma_4\sigma_3$.
It has to act as identity on $\1\otimes \1\subset V^{\otimes 3}\otimes V^{\otimes 3}$.
This is no longer the case if we replace $R_V$ by $-R_V$. This finishes the proof.

\subsection{Restriction of eigenvalues} It was shown in \cite{TW}, Section 3 that
dimensions of objects in ribbon tensor categories can be determined from braid representations
under certain circumstances. More precisely, if $Z$ is a selfdual object
in a ribbon tensor category such that 
$Z^{\otimes 2}=\bigoplus_{i=1}^k Y_i$ 
with $k\leq 5$ and the $Y_i$ mutually non-isomorphic simple objects, then
the quotient $\dim(Y_i)/\dim(Z)^2$ can be determined from the 
representation of $B_3$ on $Hom(Z,Z^{\otimes 3})$,
see the corollary in \cite{TW}, Section 3.2.   If
one of the $Y_i$s is isomorphic to $Z$, this determines the dimension of $Z$,
and hence also the dimensions of the objects $Y_i$, $1\leq i\leq k$.

\begin{theorem}\label{norootsofunity} Let $\Ca$ be a tensor category of type $G_2$
such that the image of $B_3$ in $\End_\Ca(V^{\otimes 3})$ generates
the whole algebra. Then the dimension of any object has to coincide with
the corresponding object in $\UU =Rep(U_q\g(G_2))$.
In particular, the eigenvalues of $c_{V,V}$ have to be as in
the case of the category $\UU=Rep(U_q\g(G_2))$ with $q$ not a root of unity.
\end{theorem}

$Proof.$
It follows from the discussion before this theorem for $Z$ being the object $V$
in $\Ca$ and $k=4$ that the dimensions of $V$, $V_{\Lambda_2}$ and $V_{2\Lambda_1}$
are completely determined by the 4-dimensional irreducible representation
of $B_3$ in $\End_\Ca(V^{\otimes 3})$. Hence the dimensions
of the objects $V=V_{\Lambda_1}$, $V_{\Lambda_2}$ and $V_{2\Lambda_2}$
are the same as for the quantum group $U_q(\g(G_2))$ by Proposition \ref{eigenvalues}.
As $\g(G_2)$ has rank 2, it is well-known that the dimension of any object in $\Ca$
is determined by the dimensions of the fundamental objects
 $V_{\Lambda_i}$, $i=1,2$.  This shows the first
claim. For the second claim, it follows from the explicit dimension formula
\ref{dimensionformula} that we would find a simple object $V_\la$
whose dimension would be equal to 0 for $q$ a root of unity $\neq \pm 1$. 
This would contradict rigidity of $\Ca$. 
If $q=\pm 1$, the braid representation factors through the symmetric group.
As $\dim \C S_3= 6 <\dim \End_\Ca(V^{\otimes 3})$, this is not possible under 
our assumptions.

\section{Rigidity of path representations}

\subsection{Main result of section} We consider {\it path representations of the braid groups
$B_n$ of type   $G_2$}, i.e. braid representations for the path algebra $\P$
generated by the 7-dimensional irreducible representation 
$V$ of the Lie algebra $\g(G_2)$. Recall that if $\rho_{\nu,m}$ denotes the
representation of $B_m$ on the path space $W(\nu, m)$, we have the restriction rule
\begin{equation}\label{braidrestrict}
(\rho_{\la,n})_{|B_{n-1}}\cong \bigoplus_{\mu\leftrightarrow\la} \rho_{\mu,n-1},
\end{equation}
where the summation goes over all $\mu$ for which $\la-\mu$ is a weight of $V=V_{\La_1}$,
with exceptions for the zero weight, see Section \ref{tensorsect}.
 We moreover assume the following:

\renewcommand{\labelenumi}{\alph{enumi})}
\begin{enumerate}
\item The eigenvalues of $\rho(\sigma_1)$ are as in Proposition \ref{eigenvalues} for $q$ not a root of unity.
\item The representations of $B_3$ are irreducible for all modules $W(\nu,3)$.
\item The braid $\Delta_n^2$ acts as a scalar on each module $W(\nu,n)$ compatible with a ribbon braid
structure, see Eq \ref{twistbraid} and \ref{scalarmatch}.
\end{enumerate}
We can also obtain results for path representations with less restrictive conditions, 
see Corollary \ref{scalingremark} and the remark after it. As we do not know
any non-trivial examples for this more general setting,
the goal of this section will be to prove the following theorem:

\begin{theorem}\label{maintheorem}
Any path representation of braid groups of type $G_2$ satisfying the conditions stated in
this section is uniquely determined by the eigenvalues of $\sigma_1$.  More precisely, if we
have two such path representations $\rho_1$ and $\rho_2$ such that their restrictions
to $B_2$ coincide, then also their representations of $B_n$ on any module $W(\nu,n)$
are isomorphic.
\end{theorem}

\medskip
\subsection{Outline of proof}\label{outline} 1. We will show that if we have two path representations $\rho_1$ and $\rho_2$
of type $G_2$ on modules $W_i(\nu,n)$, $i=1,2$  which satisfy the conditions stated before Theorem \ref{maintheorem}, 
then we can also achieve that the matrices for $\rho_1(\sigma_i)$ and $\rho_2(\sigma_i)$, $1\leq i<n$
coincide after rescaling the basis vectors of, say, the module $W_2(\nu,n)$. This will be done by
induction on $n$, using the restriction rule \ref{braidrestrict} as follows (with $n=2$ true by assumption):

2. To prove the claim for the module $W_2(\la, n)$, we can assume by induction assumption, using
\ref{braidrestrict} that the matrices for $\sigma_i$ with $1\leq i<n-1$ coincide for both $\rho_1$ and $\rho_2$.
Moreover these matrices  do not change if we multiply the vectors for the basis for $\rho_{\mu, n-1}$
by a non-zero scalar, say $c_\mu$, for each $\mu$.  Hence the claim will follow if we can show that
we can find suitable scalars for the basis of the $\rho_2$ representation on $W_2(\la, n)$ such that
$\rho_1(\sigma_{n-1})=\rho_2(\sigma_{n-1})$.

3.  If a matrix block for the new generator $\sigma_{n-1}$ goes through the diagrams $\mu_1,\ ...,\ \mu_r$
at level $n-1$, the matrix for $\sigma_{n-1}$ in that block will be replaced by the same  matrix conjugated by
the diagonal matrix diag$(c_{\mu_i})$ after the rescaling described in 2. So fixing this particular matrix block for $\sigma_{n-1}$
will fix the scalars $c_{\mu_i}$ in 2 (up to a common multiple). Observe that our previous results show that such a matrix block
is uniquely determined up to such a conjugation
 if the block has at most three distinct eigenvalues, see Prop. \ref{AB2representations} and Lemma \ref{smallblockss}.

4. Given two extensions $\rho_1, \rho_2$ of $ \bigoplus_{\mu\leftrightarrow\la} \rho_{\mu,n-1}$,
 we choose a block of $\sigma_{n-1}$ of
maximum size for which the matrix has at most three eigenvalues. As mentioned in 3, we can assume that the
matrices in both extensions will be the same for this block. We will then first show that we can find a  renormalization
of basis vectors such that all blocks of $\sigma_{n-1}$ with $\leq 3$ eigenvalues coincide for both extensions.
Finally, if there is a block with more than three eigenvalues, we will deduce the same result for it from the equality
of all the other blocks, using the braid relations.

\subsection{Checking braid relations} In order to check the braid relations for $\sigma_{n-2}$ and $\sigma_{n-1}$,
we consider submodules $W_3(\gamma,\lambda)$ whose basis is spanned by all paths of lengths 3 from $\gamma$ to $\lambda$.
As $\gamma$ and $\lambda$ are fixed by both  $\sigma_{n-2}$ and $\sigma_{n-1}$,  the basis vectors are 
 given by pairs $(\alpha,\beta)$ with $\alpha$ the weight on level $n-2$ and $\beta$ the weight on level $(n-1)$. 
The block $B_{\alpha}$ of $\sigma_{n-1}$ is determined by all paths with fixed first coordinate $\alpha$ 
and the block $C_{\beta}$ of $\sigma_{n-2}$
 is determined by all paths with fixed weight second coordinate $\beta$. 

\begin{remark}\label{calculationremark} In the following we will calculate a matrix entry $[\sigma_{n-1}]_{x,y}$
by exhibiting paths $s'$ and $s$ such that $[\sigma_{n-1}]_{x,y}$ is the only unknown
entry in the calculation 
$$[\sigma_{n-1}\sigma_{n-2}\sigma_{n-1}]_{s',s}=[\sigma_{n-2}\sigma_{n-1}\sigma_{n-2}]_{s',s},$$
where we only need to make sure that the other entries by which  $[\sigma_{n-1}]_{x,y}$
is multiplied are nonzero for our choice of parameters. 
\end{remark}
\begin{remark}\label{calculationremark2}
An efficient strategy for calculating the matrix
coefficients of, say, the left hand side, is as follows.
We call any sequence of paths of the form
$$s'=(\al,\beta) \quad -\quad t'=(\al,\gamma)\quad -\quad t=(\kappa,\gamma)\quad -\quad s=(\kappa,\delta)$$
a {\it 212 chain} from $(\al,\beta)$ to $ (\kappa,\delta)$.  Then it is clear that
$$[\sigma_{n-1}\sigma_{n-2}\sigma_{n-1}]_{s',s}=\sum_\gamma [\sigma_{n-1}]_{s',t'}[\sigma_{n-2}]_{t',t}[\sigma_{n-1}]_{t,s},$$
where the summation goes over all $\gamma$ which generate a $212$ chain from $s'=(\al,\beta)$ to $ s=(\kappa,\delta)$.
The calculation of $[\sigma_{n-2}\sigma_{n-1}\sigma_{n-2}]_{s',s}$ can be similarly done via $121$ chains,
where we first change the first coordinate of $s'=(\al,\beta)$; see the proof of Lemma \ref{blocks1} for an example.
\end{remark}

 The following lemma is a fairly straightforward consequence
of the braid relations. It is useful as the right hand side of the equation only includes diagonal entries of $\sigma_{n-1}$ which are
easy to calculate.

\begin{lemma}\label{blocks1}  Let $B_{\Lambda}$ and $B_{M}$ be two blocks of $\sigma_{n-1}$,
and let $s=(\Lambda,\beta)$ and $s'=(M,\beta)$. Then

$$\sum_\gamma [\sigma_{n-1}]_{s',t'}[\sigma_{n-2}]_{t',t}[\sigma_{n-1}]_{t,s}=\sum_{u\in C_\beta} [\sigma_{n-2}]_{s',u}[\sigma_{n-1}]_{u,u}[\sigma_{n-2}]_{u,s},$$
where $\gamma$ in the first sum is  such that  $t=(\Lambda,\gamma)\in B_\Lambda$ and $t'=(M,\gamma)\in B_M$.
\end{lemma}

$Proof.$ It follows from the discussion before the lemma that the left hand side is equal to
$[\sigma_{n-1}\sigma_{n-2}\sigma_{n-1}]_{s',s}$. To calculate
$[\sigma_{n-2}\sigma_{n-1}\sigma_{n-2}]_{s',s}$, just observe that the corresponding $121$ chains have to be of the form
$$s'=(\al,\beta) \quad -\quad u=(\kappa,\beta)\quad -\quad v=(\kappa,\beta)\quad -\quad s=(M,\beta);$$
hence $u=v$, which implies the claim.

\subsection{Set-up for calculating matrix entries}\label{setup} By definition
of our representations, it suffices to check the relations for $\sigma_{n-2}$ and $\sigma_{n-1}$
on subspaces $W_3(\gamma,\lambda)$ spanned by  paths of length 3 from a fixed diagram $\gamma$ to $\lambda$,
for each suitable $\gamma$. 
We will do it in detail for the most complicated case, for a weight $\lambda=a\Lambda_1+b\Lambda_2$ sufficiently
far from the walls of the Weyl chamber; this is satisfied if $a,b\geq 3$.
For other cases see Section \ref{boundarycases}.
One can check that we get the maximum number of paths for a subspace $W_3(\gamma,\lambda)$
if $\gamma=\mu$ has distance 1 from $\lambda$. 
Let us pick $\mu=(a-1)\Lambda_1+b\Lambda_2$.  Using the restriction rule for representations,
we obtain 24 paths of length 3 from $\mu$ to $\la$. These will
involve the following weights near $\lambda$ and $\mu$:
\smallskip
\vskip .4cm
\begin{center}

  \begin{tikzpicture}[scale=.4]
    \foreach \y in {5,...,7} \foreach \x in {4,...,\y}
    \draw[xshift=\x cm,thick,blue,fill=blue](\x cm,3*\x cm - 6cm) circle (.2cm);
   \foreach \y in {4,...,6} \foreach \x in {2,...,\y}
    \draw[xshift=\x cm,thick,blue,fill=blue](\x cm,3*\x cm) circle (.2cm);
       \foreach \x in {0,...,5}
    \draw[xshift=\x cm,thick,blue,fill=blue] (\x cm,3*\x cm +6 cm) circle (.2cm);
     \foreach \x in {-1,...,3}
\draw[xshift=\x cm,thick,blue,fill=blue] (\x cm,3*\x cm +12 cm) circle (.2cm);
     \foreach \x in {-2,...,1}
\draw[xshift=\x cm,thick,blue,fill=blue] (\x cm,3*\x cm +18 cm) circle (.2cm);
              \draw[magenta] (4.2 cm, 12 cm) node[anchor=north]  {\tiny $\mu$};
              \draw[magenta] (6.2 cm, 15 cm) node[anchor=north]  {\tiny $\lambda$};
              \draw (8.2 cm, 12 cm) node[anchor=north]  {\tiny $\alpha_1$};
              \draw (10.2 cm, 15 cm) node[anchor=north]  {\tiny $\alpha_2$};
              \draw (8.2 cm, 18 cm) node[anchor=north]  {\tiny $\alpha_3$};
              \draw (4.2 cm, 18 cm) node[anchor=north]  {\tiny $\alpha_4$};
              \draw (2.2 cm, 15 cm) node[anchor=north]  {\tiny $\alpha_5$};
              \draw (0.2 cm, 12 cm) node[anchor=north]  {\tiny $\alpha_6$};
              \draw (2.2 cm, 9 cm) node[anchor=north]  {\tiny $\alpha_7$};
              \draw (6.2 cm, 9 cm) node[anchor=north]  {\tiny $\alpha_8$};
  \end{tikzpicture}
\end{center}
\ignore{
$\alpha_1=(\la_1,\la_2-1)$, \\
$\alpha_2=(\la_1+1,\la_2-1)$, \\
$\alpha_3=(\la_1+1,\la_2-1)$, \\
$\alpha_4=(\la_1,\la_2+1)$, \\
$\alpha_5=(\la_1-1,\la_2+1)$, \\
$\alpha_6=(\la_1-2,\la_2+1)$, \\
$\alpha_7=(\la_1-2,\la_2)$, \\
$\alpha_8=(\la_1-2,\la_2)$.}

\begin{align}
& \alpha_1=(a+1)\Lambda_1+(b-1)\Lambda_2,  & \quad\alpha_2 =(a+2)\Lambda_1+(b-1)\Lambda_2, \cr
& \alpha_3=(a+1)\Lambda_1+b\Lambda_2, & \quad  \alpha_4 =(a-1)\Lambda_1+(b+1)\Lambda_2, \cr
& \alpha_5=(a-2)\Lambda_1+(b+1)\Lambda_2, & \quad  \alpha_6=(a-3)\Lambda_1+(b+1)\Lambda_2, \cr
& \alpha_7=(a-2)\Lambda_1+b\Lambda_2,  & \quad \alpha_8 =a\Lambda_1+(b-1)\Lambda_2.
\end{align}

\ignore{\textbf{General case (7 by 7 block):} $\lambda=a\Lambda_1+b\Lambda_2$ with $a+2b=n-1$ and $a,b\neq 0,1$. On the level $n-3$ we have $\mu=(a-1)\Lambda_1+b\Lambda_2$. \\
Let now $\alpha_1=(a+1)\Lambda_1+(b-1)\Lambda_2$, \\
$\alpha_2=(a+2)\Lambda_1+(b-1)\Lambda_2$, \\
$\alpha_3=(a+1)\Lambda_1+b\Lambda_2$, \\
$\alpha_4=(a-1)\Lambda_1+(b+1)\Lambda_2$, \\
$\alpha_5=(a-2)\Lambda_1+(b+1)\Lambda_2$, \\
$\alpha_6=(a-3)\Lambda_1+(b+1)\Lambda_2$, \\
$\alpha_7=(a-2)\Lambda_1+b\Lambda_2$, \\
$\alpha_8=a\Lambda_1+(b-1)\Lambda_2$.\\}

\medskip
\noindent
In this case $\sigma_{n-1}$ has 7 blocks given by path bases 	

$$\B_\lambda=\{t_1,t_2,t_3,t_4,t_5,t_6,t_7\},\text{ where } t_i=(\lambda, \rho) \text{ with } \rho=\lambda, \mu, \alpha_1, \alpha_2, \alpha_3, \alpha_4,\alpha_5;$$

$$\B_{\alpha_1}=\{t_8,t_9,t_{10},t_{11}\},\text{ where } t_i=(\alpha_1, \rho) \text{ with } \rho=\lambda, \mu, \alpha_1, \alpha_2;$$

$$\B_{\alpha_5}=\{t_{12},t_{13},t_{14},t_{15}\},\text{ where } t_i=(\alpha_5, \rho) \text{ with } \rho=\lambda, \mu, \alpha_4, \alpha_5;$$

$$\B_{\mu}=\{t_{16},t_{17},t_{18},t_{19}\},\text{ where } t_i=(\mu, \rho) \text{ with } \rho=\lambda, \mu, \alpha_1, \alpha_5;$$

$$\B_{\alpha_6}=\{t_{20},t_{21}\},\text{ where } t_i=(\alpha_6, \rho) \text{ with } \rho=\mu, \alpha_5;$$

$$\B_{\alpha_8}=\{t_{22},t_{23}\},\text{ where } t_i=(\alpha_8, \rho) \text{ with } \rho=\mu, \alpha_1;$$

$$\B_{\alpha_7}=\{t_{24}\},\text{ where } t_i=(\alpha_7, \rho) \text{ with } \rho=\mu.$$

\medskip
\noindent
The corresponding blocks of $\sigma_{n-2}$ are 

$$C_{\lambda}=\{t_1,t_8,t_{12},t_{16}\},$$ $$C_{\mu}=\{t_2,t_{10},t_{13}, t_{17},t_{20}, t_{22},t_{24}\},$$ $$C_{\alpha_1}=\{t_3,t_9,t_{18},t_{23}\},$$ $$C_{\alpha_2}=\{t_4,t_{11}\},$$ $$C_{\alpha_3}=\{t_5\},$$ $$C_{\alpha_4}=\{t_6,t_{14}\},$$ $$C_{\alpha_5}=\{t_7,t_{15},t_{19}, t_{21}\}.$$

\medskip
\begin{remark}\label{reflection} 
It is easy to see from the picture that we get similar block structures for the module $W_3(\al_1,\la)$,
where we just reflect the paths for each block above at the vertical axis going
through $\lambda$. Observe that the block $B_\la$ remains unchanged by this, 
as the action of $\sigma_{n-1}$ only depends on the labels $t(i)$, $n-2\leq i\leq n$
of a path $t$.
\end{remark}

\subsection{Calculating diagonal entries}
By Lemma \ref{smallblockss}, the blocks of $\sigma_{n-1}$ of size $\leq 4$ can be calculated up to conjugation by a diagonal matrix,
i.e. up to rescaling of basis vectors. This determines the diagonal entries of all such blocks as well as the product of
transposed entries, say $a_{st}a_{ts}$. We now show that we can also calculate the diagonal entries of the big 
$7\times 7$ block of $\sigma_{n-1}$.

\begin{lemma}\label{diagonalentries}
The diagonal entries of the matrices for $\sigma_{n-1}$ are uniquely determined by the entries of $\sigma_{n-2}$ and by
the results in Lemma \ref{smallblockss}.
\end{lemma}

$Proof.$ By Lemma \ref{smallblockss}, we only need to consider blocks in which $\sigma_{n-1}$ acts with more than three eigenvalues.
The only such block in $W_3(\mu,\la)$ is the block $B_\la$.
 We first find the three diagonal matrix entries $[\sigma_{n-1}]_{u,u}$ with $u=(\lambda, \beta)\in B_\lambda$ such
that $(\mu,\beta)\in B_\mu$ and $\beta\neq \mu$, i.e. for $\beta\in\{ \la, \al_1, \al_5\}$.
For this, we use Lemma \ref{blocks1} with $s=s'=(\mu,\beta)$, $\Lambda=M=\mu$,
from which we get

$$\sum_{\substack{t\in \B_{\mu} \\}} [\sigma_{n-1}]_{s,t}[\sigma_{n-2}]_{t,t}[\sigma_{n-1}]_{t,s}=\sum_{u\in C_\beta} [\sigma_{n-2}]_{s,u}[\sigma_{n-1}]_{u,u}[\sigma_{n-2}]_{u,s}.$$
The  matrix entries on the left hand side are either known diagonal entries, or products of transposed entries which are known.
The same applies to all entries  on the right hand side, except one, namely $[\sigma_{n-1}]_{u,u}$ with $u= (\lambda, \beta)\in B_\lambda$. 
 Moreover, the entry $[\sigma_{n-1}]_{u,u}$ is multiplied by matrix entries of $\sigma_{n-2}$
which are nonzero for $q$ not a root of unity, by Lemma \ref{smallblockss}.
Thus, we find $[\sigma_{n-1}]_{t_1,t_1}$,  $[\sigma_{n-1}]_{t_3,t_3}$, $[\sigma_{n-1}]_{t_7,t_7}$. 
Using the reflection symmetry, see Remark \ref{reflection},
we can similarly also calculate the diagonal entries of $\sigma_{n-1}$ for the
paths $t_1=(\la,\la)$, $t_2=(\la,\mu)$ and $t_4=(\la, \al_2)$.

To calculate the diagonal entry for the path $(\la,\al_4)$, we similarly use Lemma \ref{blocks1} for $s=s'=(\al_5,\al_4)$.
The entry for $(\la,\al_3)$ is obtained by essentially the same calculation after using the reflection symmetry in Remark \ref{reflection}.

\subsection{Determining most matrix blocks} We are dealing with the first part of point 4 of our outline. This means we are
going to show that there exists a diagonal matrix which conjugates each block of $\rho_1(\sigma_{n-1})$
in which the matrix acts with $\leq 3$ eigenvalues to the corresponding block of $\rho_2(\sigma_{n-1})$
and does not change the image of $B_{n-1}$.

\begin{lemma}\label{smallblocks} Let $\rho_1$ and $\rho_2$ be two representations of $B_n$ with the same path basis such that 
${\rho_1}_{|_{B_{n-1}}}={\rho_2}_{|_{B_{n-1}}}$ produce the same matrices. 
Then we can also make coincide all  blocks for $\sigma_{n-1}$ in which it has at most three different eigenvalues.
\end{lemma} 

$Proof.$ By point three of our outline, we can assume that the matrices for block $B_\mu$
coincide for both representations. We will show that the braid relations will essentially determine the action
of the generator $\sigma_{n-1}$ for the blocks in the statement. We will frequently
use the following two observations:
\begin{enumerate}
\item  The right hand side of Lemma \ref{blocks1} can always be calculated,  by Lemma \ref{diagonalentries}.
\item  If the matrix entry $[\sigma_{n-1}]_{x,y}$ can be calculated from the matrix entry 
$[\sigma_{n-1}\sigma_{n-2}\sigma_{n-1}]_{s',s}$ as outlined in Remark \ref{calculationremark},
we can similarly calculate the entry  $[\sigma_{n-1}]_{y,x}$   from the matrix entry 
$[\sigma_{n-1}\sigma_{n-2}\sigma_{n-1}]_{s,s'}$.
\end{enumerate}
$Notation$: We will say that a matrix entry {\it is known} if it can be expressed in terms of entries
of $\sigma_{n-2}$ and of the entries of $\sigma_{n-1}$ from the fixed block $B_\mu$.
Observe that this means that such entries have to coincide in the two path representations $\rho_1$ and $\rho_2$.

\begin{enumerate}
\item We first find the matrix entries for the two $2\times 2$ blocks $B_\alpha$ with $\alpha=\alpha_6$ or $\alpha_8$. 
We outline the calculation for $\al=\al_6$, where we use Lemma \ref{blocks1} with $s=(\al_6, \mu)$ and $s'=(\mu, \mu)$. 
The only unknown quantity on the left hand side is  $[\sigma_{n-1}]_{t,s}$
for $t=(\al_6,\al_5)$; the only other possibility for $t$ would be $t=s$ for which we get a known
diagonal entry  $[\sigma_{n-1}]_{s,s}$. Hence we can solve for $[\sigma_{n-1}]_{t,s}$, by Lemma \ref{smallblockss}.

\item Let now $B_\alpha$ be one of the two other $4 \times 4$ blocks, i.e. $\alpha=\alpha_1$ or $\alpha_5$. 
We do the case with $\al = \al_1$. We first calculate the entry  $[\sigma_{n-1}]_{t,s}$, where $t= (\alpha_1, \mu)$
and $s= (\alpha_1, \alpha_1)$ by using Lemma \ref{blocks1} with $s'=(\al_8,\mu)$. As $B_{\al_8}$
only has two paths, we only have two summands on the left hand side in that Lemma, with the only unknown quantity
$[\sigma_{n-1}]_{t,s}$.

Next, we consider matrix entries involving the path $r=(\alpha_1, \lambda)$. We again use Lemma \ref{blocks1} with $s=(\alpha_1,\rho)$, $s'=(\mu, \rho)$, where $\rho\in\{ \alpha_1, \mu\}$.
There are three summands on the left hand side, for which all quantities are known except $[\sigma_{n-1}]_{r,s}$
(the other two quantities $[\sigma_{n-1}]_{t,s}$ are either a diagonal entry or known from the previous paragraph).

Recall that all matrix blocks where the braid generator has $\leq 3$ eigenvalues are already determined by Lemma \ref{AB2representations}, up to conjugation
by a diagonal matrix. Now we only need to observe that the  path $t_{11}=(\al_1,\al_2)$
belongs to a block of $\sigma_{n-2}$ which does not have paths appearing in $B_\mu$.  Hence we can rescale the paths in the block $C_{\al_2}$
so that the block $B_{\alpha_1}$ is the same in both representations $\rho_1$ and $\rho_2$.
\end{enumerate}

\subsection{Calculating the unknown block} We have seen in Lemma \ref{smallblocks} that we can make equal all but one matrix block
of $\rho_1(\sigma_{n-1})$ with  $\rho_2(\sigma_{n-1})$ on $W_3(\mu,\lambda)$  if all the blocks of $\rho_1(\sigma_{n-2})$ 
equal the ones of $\rho_2(\sigma_{n-2})$ and there is one  $4\times 4$ block for which  $\rho_1(\sigma_{n-1})$ equals  $\rho_2(\sigma_{n-1})$. Moreover, we also know that the diagonal entries have to coincide for all blocks.
We now want to show that this also enables us to make the remaining matrix entries equal.

\ignore{
\item We now find the $4 \times 4$ subblock of the $7 \times 7$ block $B_\lambda$ given by paths $(\lambda,\beta)$ for which $(\mu,\beta)$ is a path in $B_\mu$. 
In this case, we have $\beta\in\{ \la, \mu, \alpha_1, \alpha_5\}$.
Pick two such paths $s'_i=(\lambda, \beta_i)$, $i=1,2$ and let $s_i=(\mu, \beta_i)$. 
To find $[\sigma_{n-1}]_{s'_1, s'_2}$, we look at
$$\sum_{\substack{u\in \B_{\mu} }} [\sigma_{n-1}]_{s_1,u}[\sigma_{n-2}]_{u,u}[\sigma_{n-1}]_{u,s_2}=\sum_{\substack{v=(\gamma,\beta_1)\in C_{\beta_1},\\ \text{such that}\\w=(\gamma, \beta_2)\in C_{\beta_2}} }[\sigma_{n-2}]_{s_1,v}[\sigma_{n-1}]_{v,w}[\sigma_{n-2}]_{w,s_2}.$$
We can calculate the left hand side as all matrix entries $[\sigma_{n-1}]_{s_1,u}$ and $[\sigma_{n-1}]_{u, s_2}$ are from the block $B_\mu$. 
As we already know the diagonal entries of $\sigma_{n-1}$ in this particular subblock, we can assume that not both paths $s'_1$ and $s'_2$ are of the form
$(\lambda,\mu)$. Hence we have at most four summands on the right hand side. We can then check that 
$v=(\beta, )$ is the only matrix entry from $B_\lambda$ is $[\sigma_{n-1}]_{s_1',s_2'}$. 
\item Next we find the quotients in the remaining entries of $4 \times 4$ blocks $B_\alpha$ with $\alpha=\alpha_1$ or $\alpha_5$ using Lemma \ref{blocks2}. 
For $\alpha=\alpha_1$, let $a=t_{11}$ and for $\alpha=\alpha_5$, let $a=t_{14}$. Then we can find $[\sigma_{n-1}]_{p,a}[\sigma_{n-1}]_{a,s}$ from Lemma \ref{blocks2} using part (d). We will let $\Lambda=\mu$, and $p,s\in B_\alpha$ with $p,s\neq a$. This gives all the quotients $$\frac{[\sigma_{n-1}]_{p,a}}{[\sigma_{n-1}]_{s,a}}\text{ and } \frac{[\sigma_{n-1}]_{a,p}}{[\sigma_{n-1}]_{a,s}}.$$
\item Now we can find the matrix entries $[\sigma_{n-1}]_{p,r}$ for $p$ or $r=t_4$ or $t_6$ using the same method above and Lemma \ref{blocks1} depending on the parameter of quotients in part (e). This is because $t_4'=t_{11}$ in $B_{\alpha_1}$ and $t_6'=t_{14}$ in $B_{\alpha_5}$.
\item To find the remaining matrix entries $[\sigma_{n-1}]_{p,r}$ for $p$ or $r=t_7$, we use Lemma \ref{blocks2} again to find 
$[\sigma_{n-1}]_{t,a}[\sigma_{n-1}]_{a,s}$ with $a=t_7$. Then we can find all the quotients: $$\frac{[\sigma_{n-1}]_{p,t_7}}{[\sigma_{n-1}]_{s,t_7}}\text{ and } \frac{[\sigma_{n-1}]_{t_7,p}}{[\sigma_{n-1}]_{t_7,s}}.$$}

\begin{lemma}\label{bigblock} Let $s_1$ and $s_2$ be paths  belonging to the unknown
$7 \times 7$  block $B_\la$ of $\sigma_{n-1}$. Then the matrix entry $[\sigma_{n-1}]_{s_1, s_2}$ is uniquely determined
if none of the paths $s_1$ or $s_2$  is of the form $(\la,\mu)$ or $(\la,  \al_3)$.
\end{lemma}
$Proof.$ As neither of the paths $s_1$, $s_2$ belongs to the $1\times 1$ block $C_{\al_3}$, we can find 
paths $t_i\neq s_i$ in the same $\sigma_{n-2}$ block. As neither path $t_i$ belongs to the unknown block $B_\la$
we can calculate
$$[\sigma_{n-1}\sigma_{n-2}\sigma_{n-1}]_{t_2,t_1}=\sum_{u,v} [\sigma_{n-1}]_{t_2,u}[\sigma_{n-2}]_{u,v}[\sigma_{n-1}]_{v,t_1}.$$
As neither path $s_i$ or $t_i$ is in the block $C_{\mu}$ by assumption,
we also know that the entries  $[\sigma_{n-2}]_{t_i, s_i}\neq 0$ for $i=1,2$, by Lemma \ref{smallblockss}.
Hence we can solve for $[\sigma_{n-1}]_{s_2, s_1}$, provided we can show it is the only unknown
entry in
$$[\sigma_{n-2}\sigma_{n-1}\sigma_{n-2}]_{t_2,t_1}=\sum_{u',v'} [\sigma_{n-2}]_{t_2,u'}[\sigma_{n-1}]_{u',v'}[\sigma_{n-2}]_{v',t_1}.$$
Writing $\sigma_{n-2}\cdot t_i$ as a linear combination of paths $t$, there is only one of them
in $B_\lambda$, namely $s_i$. Hence all other matrix coefficients $[\sigma_{n-1}]_{u',v'}$ in 
the sum above belong to blocks whose coefficients are known.
This finishes the proof.

\ignore{
\begin{theorem}\label{braidunique} Any braid representations on $G_2$ paths is
uniquely determined by the eigenvalues of $\sigma_1$ and the choice of a
square root of the determinant of $\sigma_1$ in the 4-dimensional
irreducible representation of $B_3$.
\end{theorem}}

\medskip

In order to calculate the additional matrix entries of $B_\la$, we will use Remark \ref{reflection}.
Applying Lemma \ref{smallblocks} to all paths in $W_3(\al_1,\la)$ we can similarly
assume that all corresponding blocks $B_\gamma$ of $\rho_i(\sigma_{n-1})$ of size $\leq 4$ also coincide on $W_3(\al_1,\la)$
for $i=1,2$. Moreover, we can also assume that if a block $B_\gamma$ appears in both $W_3(\mu,\la)$ and $W_3(\al_1,\la)$,
we get the same matrices.

\begin{lemma}\label{bigblock2}  After suitable renormalizations of the path basis vectors for $\rho_2$,
the matrices $\rho_1(\sigma_{n-1})$ and $\rho_2(\sigma_{n-1})$ coincide.
\end{lemma}

$Proof.$  
By Lemma \ref{smallblocks}, it only reminds to check the claim for those matrix entries in $B_\la$ which have
not already been covered in Lemma \ref{bigblock}. Observe that if we apply Lemma \ref{bigblock} to the block $B_\la$
in $W_3(\al_1,\la)$, we have to exclude the paths $(\la,\al_1)$ and $(\la, \al_4)$, by using the symmetry
in Remark \ref{reflection}. Hence we can calculate all matrix entries in $B_\la$ except those entries 
$[\sigma_{n-1}]_{u,v}$ where one path is  $t_2=(\la,\mu)$ or $t_5=(\la,\al_3)$, and the other path
is   $t_3=(\la,\al_1)$ or  $t_6=(\la, \al_4)$. There are eight such entries, which always come in pairs
$(u,v)$ and $(v,u)$. We proceed as outlined in Remark \ref{calculationremark}. As observed in the proof 
of Lemma \ref{smallblocks}, it suffices to indicate how to calculate one entry for each of these pairs.
Here we use notation from Section \ref{setup}:

For $(u,v)=(t_2,t_3)$, we take $(s',s)=(t_2, t_{22})$;\\
For $(u,v)=(t_2,t_6)$, we take $(s',s)=(t_2, t_{13})$;\\
For $(u,v)=(t_5,t_3)$, we take $(s',s)=(t_5, t_{10})$;\\
For $(u,v)=(t_5,t_6)$, we take $(s',s)=(t_5, t_{13})$.

\subsection{Matrices for weights near the boundary of the Weyl chamber}\label{boundarycases}
 If the weight $\la$ is near the
boundary of the Weyl chamber, we can calculate the block $B_\la$ similarly as before. However, some of
the arguments need to be changed. In particular, we can usually not apply Remark \ref{reflection}.
The most involved case is if we take $\la= a\Lambda_1$ and $\mu=(a-1)\Lambda_1$. We can still use
the picture for the weights in the general case, except we need to remove all weights below
the line connecting $\la$ and $\mu$. One can see from this that the block $B_\la$ has 5 paths.
We proceed as in the general case, where we need to make the following changes (the necessary
lemmas appear in the next section).

\begin{enumerate}
\item We can calculate three diagonal entries of the block $B_\la$ as in Lemma \ref{diagonalentries}, 
and we can calculate the diagonal entry for
the path $(\la,\mu)$ from Lemma \ref{projection}. This also determines the fifth diagonal 
entry as we know the trace of $\sigma_{n-1}$ in that block.

\item We can  then calculate all matrix blocks of size $\leq 4$ as in Lemma \ref{smallblocks},
assuming we fix the matrix for block $B_{\al_5}$.

\item We calculate all matrix entries of the block $B_\la$ not involving the paths $(\la, \mu)$ and $(\la, \al_3)$
as in Lemma \ref{bigblock}.

\item We calculate the matrix entries $[\sigma_{n-1}]_{u,v}$  involving the path $u=(\la, \mu)$ 
but not the path $(\la, \al_3)$
similarly as it was done at the end of the proof of Lemma \ref{bigblock2} as follows
(it is important to do it in the given order):

For $v=(\la,\al_5)$, we take $(s',s)=(u, (\al_6,\mu))$;\\
For $v=(\la,\la)$, we take $(s',s)=(u, (\mu,\mu))$;\\
For $v=(\la,\al_4)$, we take $(s',s)=(u, (\al_5,\mu))$.

\item The only missing matrix entries of the block $B_\la$ are the ones involving the
path $t_5=(\la, \al_3)$. These are determined in Lemma \ref{bigblock3}, where $t_0=t_5$.
\end{enumerate}
The blocks $B_\la$ near the left boundary can be similarly calculated, where one has to make
far fewer adjustments. E.g. if $\la = a\Lambda_2$ and $\mu = \Lambda_1 + (a-1)\Lambda_2$,
the block $B_\la$ only has three paths and can be calculated as in Lemma \ref{smallblockss}.
If $\la= \Lambda_1 + a\Lambda_2$ and $\mu = a\Lambda_2$, we can again use Remark \ref{reflection}.

\subsection{Additional lemmas} Let $p_i$ be the eigenprojection of the image of $\sigma_i$
corresponding to the eigenvalue $\al_{(0,0)}$. We see from the Bratteli diagram
that $p_1\rho(\sigma_2)^mp_1=\gamma_m p_1$   for all $m\in \Z$,
where $\gamma_m$ is a scalar.
The same relation holds if we shift the indices to $n-2$ and $n-1$.
The projection $p_{n-2}$ only is nonzero in the block $C_\mu$. For the path algebra
for $U_q\g(G_2)$ (in fact, for any self-dual representation of a quantum group)
its diagonal entries
can be calculated as it was done in e.g. \cite{RW} (2.15) in connection with orthogonal and symplectic
groups. Indeed, the diagonal entry of $p_{n-2}$ for the path $t=(\gamma,\mu)$ is given by
$$\frac{\dim_q V_\gamma}{\dim_q V\ \dim_q V_\mu},$$
with $q$-dimensions as given in Eq \ref{dimensionformula}. In particular, we also obtain from
this that the diagonal entries are non-zero for $q$ not a root of unity.
As $p_{n-2}$ acts as a rank 1 idempotent in that block, all its matrix entries in that block
are non-zero.

\begin{lemma}\label{projection} Assume we have a path representation of $B_n$
such that its restriction to $B_{n-1}$ is isomorphic to the one of $U_q\g(G_2)$.
Then the diagonal entry of $\sigma_{n-1}$ for the path $(\la,\mu)$ can be calculated from
the other diagonal entries of $\sigma_{n-1}$ and the entries of $p_{n-2}$.
\end{lemma}

$Proof$. We compare the $u=(\la,\mu)$ diagonal entries of the equation
$p_{n-2}\rho(\sigma_{n-1})p_{n-2}=\gamma_1  p_{n-2}$ to obtain the equation
$$ \sum_{\nu} p_{u, (\nu,\mu)} [\sigma_{n-1}]_{(\nu,\mu),(\nu,\mu)}  p_{ (\nu,\mu),u}\ =\ \gamma p_{u,u},$$
where $p_{u,v}$ is a matrix entry of $p_{n-2}$.
Observe that the product of opposite off-diagonal entries of $p_{n-2}$ is equal to the product of the
corresponding diagonal entries, which are known and nonzero by assumption. Hence  we can solve for the only unknown diagonal entry of $\sigma_{n-1}$
with $\nu=\la$ in that equation.

\begin{lemma}\label{bigblock3}  Assume all matrix entries of  $\rho_1(\sigma_{n-1})$ and $\rho_2(\sigma_{n-1})$
coincide except possibly the non-diagonal entries involving the path $t_0$ belonging to the $1\times 1$ block $C_{\al_3}$.
Then we can renormalize the basis vector  $t_0$ 
such that the matrices $\rho_1(\sigma_{n-1})$ and $\rho_2(\sigma_{n-1})$ coincide.
\end{lemma}
$Proof.$  
Let $t_1$ and $t_2$ be paths in $B_\la$ such that $t_i\neq t_0$ for both $i= 1,2$. 
Then we can calculate
$$[\sigma_{n-2}\sigma_{n-1}\sigma_{n-2}]_{t_1,t_2}=\sum_{u,v\neq t_0} [\sigma_{n-2}]_{t_1,u}[\sigma_{n-1}]_{u,v}[\sigma_{n-2}]_{v,t_2};$$
indeed as $u(n-1)=t_1(n-1)\neq t_0(n-1)$ and $v(n-1)=t_2(n-1)\neq t_0(n-1)$,
 $[\sigma_{n-1}]_{u,v}$ is known for all summands on the right hand side,
and hence so is the whole sum.
By Lemma \ref{blocks1}
$$[\sigma_{n-1}\sigma_{n-2}\sigma_{n-1}]_{t_1,t_2}=\sum_{s\ {\rm in}\  B_\la}[\sigma_{n-1}]_{t_1,s}[\sigma_{n-2}]_{s,s}[\sigma_{n-1}]_{s,t_2},$$
where we know all summands on the right hand side except the one with $s=t_0$. The $1\times 1$ block entry
 $[\sigma_{n-2}]_{t_0,t_0}\neq 0$. Hence we can solve for 
$[\sigma_{n-1}]_{t_1,t_0}[\sigma_{n-1}]_{t_0,t_2}$, for any $t_1, t_2$ in the unknown block
in terms of known matrix entries. Pick a $v=t_2\neq t_0$ such that $[\sigma_{n-1}]_{v,t_0}[\sigma_{n-1}]_{t_0,v}\neq 0$.
Then also the  quotients $[\sigma_{n-1}]_{t_1,t_0}/[\sigma_{n-1}]_{v,t_0}$ are determined in terms of known matrix entries.

We can now rescale the basis vector $t_0$ for $\rho_2$ such that  $\rho_1(\sigma_{n-1})_{v,t_0}=\rho_2(\sigma_{n-1})_{v,t_0}$.
Then the claim follows from the statements of the results in the previous paragraph.

\subsection{Proof\ of\ Theorem\ \ref{maintheorem}} 
We will prove the theorem for representations of $B_n$
on paths of length $n$ by induction on $n$. The statement is true by assumption for $n=2$.
For $n=3$, we have irreducible representations  of $B_3$ up to dimension 4.
Here the statement follows from \cite{TW}: It is proven there that irreducible representations
of $B_3$ up to dimension 3 are completely determined by the eigenvalues of $\sigma_1$.
For dimension 4, there are two possibilities for given eigenvalues of $\sigma_1$, depending
on the choice of the square root of the determinant, as
the central element $\Delta_3^2$ acts via the scalar $det(\sigma_1)^{3/2}$. 
This square root is determined by the ribbon braid structure
and has to be equal to $q^{-8}$, see the proof of Proposition \ref{eigenvalues}.

For the general induction step, we proceed as sketched in the outline. 
By induction assumption, we can assume that the matrices $\rho_1(\sigma_i)$
and $\rho_2(\sigma_i)$
coincide for $i<n-1$. By Lemmas \ref{smallblocks} and \ref{bigblock}
we can also achieve that the matrices $\rho_j(\sigma_{n-1})$ coincide for $j=1,2$.

\begin{corollary}\label{scalingremark}
The uniqueness result of Theorem \ref{maintheorem} would also hold for any path representation
of type $G_2$ for which the statement in Lemma \ref{smallblockss} about certain
matrix entries being nonzero holds.
\end{corollary}

$ Proof.$ The proof of Theorem \ref{maintheorem} only needed that certain matrix entries
were well-defined and nonzero. Their explicit values were irrelevant.

\begin{remark}
It is easy to obtain additional path representations from a given one
by multiplying the matrices for the standard generators $\sigma_i$
by a common constant. We are not aware of any other path representations
of type $G_2$ besides the ones already mentioned.
\end{remark}

\subsection{Application to tensor categories} Let $\Ca$ be a semisimple ribbon
tensor category whose fusion rules are the ones of $\g(G_2)$, and let
$\UU$ be equal to $Rep(U_q\g(G_2))$. The following theorem has already appeared
before in \cite{MPS}, where it was proved by a different method (see also the remarks below).

\begin{theorem}\label{tensorequivalence} Let $V$ be the object in $\Ca$
corresponding to the 7-dimensional representation of $U_q\g(G_2)$.
Assume that the image of $B_3$ generates $\End_\Ca(V^{\otimes 3})$.
Then  $\Ca\cong Rep(U_q\g(G_2))$ for some $q$ not a root of unity.
\end{theorem}

$Proof.$ It follows from Theorem \ref{norootsofunity} that the eigenvalues
of the braiding element $c_{V,V}$ coincide with the ones for the corresponding
braiding morphism in $\UU=Rep(U_q\g(G_2))$ for $q$ not a root of unity.
Theorem \ref{maintheorem} and Theorem \ref{firstfund} now imply that $\End_\Ca(V^{\otimes n})
\cong \End_\UU(V^{\otimes n})$ such that
 the tensor embeddings
$$\End_\Ca(V^{\otimes n})\otimes\End_\Ca(V^{\otimes m}) \  \to\ \End_\Ca(V^{\otimes n+m})$$
are compatible with these isomorphisms. As explained in \cite{KW}, see also \cite{TW2},
this almost implies the equivalence of categories $\Ca$ and $\UU$. The only additional
necessary information comes from the embeddings of the trivial object $\1$ into some tensor power
$V^{\otimes N}$. If we have
$$\iota_N : \1 \to V^{\otimes N}\quad {\rm and} \quad \pi_N: V^{\otimes N}\to \1 $$
such that $\pi_N \circ\iota_N $ is the identity of $\1$, we obtain a number $\Theta(N)$
defined by
$$\Theta(N)1_V\ =\ (\pi_N\otimes 1_V)\ c_{V,V^{\otimes N}}\ (1_V\otimes \iota_N).$$
One easily derives from this also that 
$$c_{V,V^{\otimes N}}\ (1_V\otimes \iota_N)\ =\ \Theta(N)(\iota_N\otimes 1_V),$$
$$c_{V^{\otimes N},V}(\iota_N\otimes 1_V)\ =\ \Theta(N)^{-1}(1_V\otimes \iota_N),$$
see \cite{TW2}, Lemma 4.1 (where part (b) is correct if $\Theta$ is replaced by $\Theta^{-1}$).
It follows from \cite{TW2}, Proposition 4.5 that $\Theta(N)^N=1$.
Now observe that the trivial representation $\1$ appears both in the second and the third
tensor power with multiplicity 1. Using the formulas above and the braiding axiom \ref{braid1}, we calculate
the middle expression in the formula below in two different ways to obtain
$$\Theta(3)^2 1_{\1}\ =\ (\pi_3\otimes \pi_2)\ c_{V^{\otimes 2}, V^{\otimes 3}}\ (\iota_2\otimes \iota_3)
\ =\ \Theta(2)^{-3}  1_{\1}.$$
Hence it follows that $\Theta(3)=1=\Theta(2)$. The equivalence of $\Ca$ and $\UU$ now follows e.g.
from \cite{TW2}, Theorem 4.8 and Corollary 3.6 or \cite{KW}, Proposition 2.1.

\begin{remark} 1. The techniques in this paper should also work for fusion tensor categories
of type $G_2$, for $q$ a root of unity, assuming that $B_3$ generates
$\End(V^{\otimes 3})$. Our approach might even be useful to prove analogous results
for tensor products of type $E_6$ and $E_7$. Indeed, results similar to the ones
in Lemma \ref{smallblockss}
have already been proved in \cite{Wexc} for Lie types $E_N$ in general. 

2. Our methods can not be applied for symmetric tensor categories. However, 
it would seem
reasonable to expect that analogous techniques, using
 infinitesimal braid relations and Casimirs, see e.g. \cite{LZ} might still work
for such cases.

3. We have results which would replace our surjectivity assumption
for the third tensor power by somewhat weaker assumptions
(e.g. just assuming semisimplicity of the representation of $B_3$
for $q\neq \pm 1$). We would obtain a complete classification of
tensor categories of type $G_2$ for $q\neq \pm 1$ if we could show that the
surjectivity assumptions would always hold in these cases.

4. Theorem \ref{tensorequivalence} has already been proved in
\cite{MPS} with the same assumptions for
$\End(V^{\otimes 3})$ as in our paper, including a generalization 
to $q=\pm 1$. The authors used
Kuperberg's spiders in \cite{Ku2} which do not play a role in our
approach.
\end{remark}

\end{document}